\documentclass[a4paper,reqno]{amsart}

\usepackage{amsmath}
\usepackage{amsthm}
\usepackage{amsfonts}
\usepackage{amssymb}

\usepackage{color}

\usepackage{hyperref}

\hypersetup{
    colorlinks=false,
    pdfborder={0 0 0},
    pdfstartview={XYZ null null 1.00}
}

\numberwithin{equation}{section}

  \def\<{\langle}
  \def\>{\rangle}

\theoremstyle{plain}
  \newtheorem{theorem}{Theorem}[section]
  \newtheorem{proposition}[theorem]{Proposition}
  \newtheorem{lemma}[theorem]{Lemma}
  \newtheorem{corollary}[theorem]{Corollary}

\theoremstyle{definition}
  \newtheorem{definition}[theorem]{Definition}
  
  \newtheorem{remark}[theorem]{Remark}

\begin{document}

\title[Effect of resonance on the existence of periodic solutions...]{Effect of resonance on the existence of periodic solutions for strongly damped wave equation}

\author{Piotr Kokocki}
\address{\noindent BCAM - Basque Center for Applied Mathematics \newline Alameda de Mazarredo 14, 48009, Bilbao, Spain}
\address{\noindent  Faculty of Mathematics and Computer Science \newline Nicolaus Copernicus University \newline  Chopina 12/18, 87-100 Toru\'n, Poland}
\email{pkokocki@mat.umk.pl}
\thanks{The researches supported by the NCN Grant no. NCN  2013/09/B/ST1/01963}



\begin{abstract}
We are interested in the differential equation
$\ddot u(t) = -A u(t) - c A \dot u(t) + \lambda u(t) + F(t,u(t))$, where $c > 0$ is a damping factor, $A$ is a sectorial operator and $F$ is a continuous map. We consider the situation where the equation is at resonance at infinity, which means that $\lambda$ is an eigenvalue of $A$ and $F$ is a bounded map. We introduce new geometrical conditions for the nonlinearity $F$ and use topological degree methods to find $T$-periodic solutions for this equation as fixed points of Poincar\'e operator.
\end{abstract}

\maketitle

\setcounter{tocdepth}{2}

\section{Introduction}

In this paper we are interested in the following strongly damped wave equation
\begin{equation}\label{eq2wave}
\left\{\begin{aligned}
    & u_{tt} -c\Delta u_t = \Delta u  + \lambda u + f(t,x,u), && t\ge 0, \ x\in\Omega \\
    & u(t,x) = 0, && t\ge 0, \ x\in\partial\Omega
\end{aligned}\right.
\end{equation}
where $c> 0$ are damping factors, $\lambda$ is a real number and $f:[0,+\infty)\times\Omega\times\mathbb{R}\to\mathbb{R}$ is a continuous map on an open bounded set $\Omega\subset\mathbb{R}^n$, which is $T$-periodic in time.

The existence of periodic solutions for damped wave equations has been investigated by many authors in the last years. In particular a large part of these studies concerns the weakly damped wave equation
\begin{equation}\label{w-damp}
\left\{\begin{aligned}
    & u_{tt} -c u_t = \Delta u  + \lambda u + f(t,x,u), && t\ge 0, \ x\in\Omega \\
    & u(t,x) = 0, && t\ge 0, \ x\in\partial\Omega
\end{aligned}\right.
\end{equation}
where the Laplacian in the damping term appears in the zero fractional power. For instance the results obtained in the series of papers \cite{MR2379473}, \cite{MR1040261}, \cite{MR1473860}, \cite{MR1646610}, \cite{MR1735817} provides the existence of $T$-periodic solutions for \eqref{w-damp}, in the case when $\Omega\subset\mathbb{R}^n$ is {\em a thin domain}, that is, a cartesian product of an open bounded subset of $\mathbb{R}^{n-1}$ and a small open interval. In these papers the periodic solutions are obtained as fixed points of the Poincar\'e operator by topological degree methods.
On the other hand, we refer the reader to \cite{MR0692283} where a homotopy invariants method is used to study the existence of periodic solutions in the case where $\Omega$ is an open interval and the damping term $u_t$ is additionally involved with a nonlinearity. See also \cite{MR2401532}, \cite{MR0847771}, \cite{MR0512156}, \cite{MR0460914}, \cite{MR1621766} for the results where $\Omega$ is again an open interval with the difference that the zero Dirichlet boundary conditions in the equation \eqref{w-damp} are replaced by the periodic one.

From the point of view of the mathematics, physics and engineering it is of importance to consider the equation \eqref{eq2wave} in the presence of {\em resonance at infinity}, which means that
$$\mathrm{Ker}\,(\lambda I - A_2)\neq \{0\} \text{ \ and  \ } f \text{  \ is a bounded map},$$
where we define $A_2 u := -\Delta u$ for $u\in D(A_2) := H^2(\Omega)\cap H^1_0(\Omega)$. Here we refer the reader to \cite{MR1084570}, \cite{MR1389086} for an extensive discussion on the meaning of resonance in the periodic oscillations of suspension bridges.
The existence of periodic solutions for the equation \eqref{eq2wave} in the case of the resonance at infinity was considered in \cite{MR0513090} under the assumption that the damping constant $c = 0$. There was proved that the equation admits a periodic solution provided the nonlinearity $f$ satisfies so called Landesman-Lazer type conditions. Subsequently, these conditions and topological degree methods were used in \cite{Aleks} to obtain the existence of periodic solutions in the weakly damped case \eqref{w-damp}.

In this paper our aim is to study the existence of $T$-periodic solutions for {\em the strongly damped wave equation} \eqref{eq2wave} in the presence of the resonance at infinity, which seems to be not explored problem so far. Throughout the paper we will consider the more general abstract differential equation
\begin{equation}\label{row-par}
\ddot u(t) = -A u(t) - c A \dot u(t) + \lambda u(t) + F(t,u(t)), \quad t \in [0,+\infty)\\
\end{equation}
where $c> 0$ is still a damping constant, $\lambda$ is a real number, $A: X\supset D(A)\to X$ is a positive sectorial operator with compact resolvents on a Banach space $X$ and $F:[0,+\infty)\times X^\alpha\to X$ is a continuous map, where $X^\alpha = D(A^\alpha)$ for $\alpha\in(0,1)$, is a fractional space endowed with the graph norm. For more details on the construction and properties of the fractional spaces we refer the reader to \cite{MR610244}, \cite{MR0089373}, \cite{Pazy}.

After passing into the abstract framework, we will say that the equation \eqref{row-par} is at {\em the resonance at infinity}, provided
$$\mathrm{Ker}\,(\lambda I - A)\neq \{0\} \text{ \ and \ } F \text{ \ is a bounded map.}$$
The main difficulty lies in the fact that, in the presence of resonance, there are examples of the nonlinearity $F$ such that the equation \eqref{row-par} do not admit a periodic solution. This fact will be explained in Remark \ref{rem-non-ex4}. To overcome this difficulty we address the naturally raising question which says: \\[-5pt]
\begin{equation}\label{aasasa}
\left.\begin{aligned}
 \parbox[t][][t]{105mm}{\em what additional assumptions for the nonlinearity $F$ should be made to prove that the equation \eqref{row-par} admits a $T$-periodic mild solution.}
\end{aligned} \ \ \right\} \\[4pt]
\end{equation}
To explain our methods more precisely, observe that the equation \eqref{row-par} can be written in the following form
\begin{equation*}
    \dot w(t) = -{\bf A} w(t) + {\bf F}(t, w(t)), \qquad t > 0,
\end{equation*}
where ${\bf A}:{\bf E}\supset D({\bf A})\to{\bf E}$ is a linear operator on the space ${\bf E}:=X^\alpha\times X$ given by
\begin{equation*}
\begin{aligned}
D({\bf A}) & :=\{(x,y)\in {\bf E} \ | \ x + c y\in D(A)\} \\
{\bf A}(x,y) & :=(-y,A(x + c y) - \lambda x) \qquad \mathrm{for} \quad (x,y)\in D({\bf A}),
\end{aligned}
\end{equation*}
and ${\bf F}:[0,+\infty)\times {\bf E}\to{\bf E}$ is a map defined by
$${\bf F}(t, (x,y)):=(0,F(t, x)) \qquad \mathrm{for} \quad t\in[0,+\infty), \ \ (x,y)\in{\bf E}$$
Let us assume that, for every initial data $(x,y)\in {\bf E}$, the equation \eqref{row-par} admits a mild solution $w(\,\cdot\,;(x,y)):[0,+\infty)\to {\bf E}$ starting at $(x,y)$. Then the $T$-periodic solutions of \eqref{row-par} can be identified with fixed points of \emph{the Poincar\'e operator} ${\bf \Phi}_T:{\bf E} \to {\bf E}$, defined by ${\bf \Phi}_T (x,y):= w(T;(x,y))$ for $(x,y)\in {\bf E}$.

Using the fact that $A$ has compact resolvent, we prove in Section 3 that the space ${\bf E}$ can be endowed with a norm such that the operator ${\bf \Phi}_T$ is condensing in the Hausdorff measure of noncompactness. Therefore, the natural way to obtain the existence of fixed points for the Poincar\'e operator is to exploit the Sadovski version of the Schauder fixed point theorem (see e.g. \cite{Dugundji-Granas}). Unfortunately, this theorem can not be used in the explicit way because it is difficult to find a bounded and convex set that is invariant under ${\bf \Phi}_T$.
Therefore we will apply the approach based on topological degree theory for condensing fields (see e.g. \cite{MR0246285}, \cite{MR0306986}, \cite{MR0312341}) that allows us to obtain effective methods to search for the fixed points of ${\bf \Phi}_T$. More precisely, to address the question \eqref{aasasa}, we introduce resonant conditions for the nonlinearity $F$ and we use them to construct an open bounded set $N$ with the property that its boundary does not contain fixed points of the Poincar\'e operator. Then we compute the topological degree of $I - {\bf \Phi}_T$ with respect to $N$, in the terms of the resonant conditions imposed earlier on the nonlinearity. The non-triviality of this degree will guarantee the existence of $T$-periodic solutions.

Roughly speaking, the resonant conditions read as follows. Let $\lambda\in\mathbb{R}$ be given eigenvalue of the operator $A$ such that its the geometric and algebraic multiplicities are equal. Since $A$ has compact resolvents, the eigenvalue $\lambda$ is isolated and there is a subspace $V\subset X$ with the property that $X = \mathrm{Ker}\,(\lambda I - A)\oplus V$ and $\sigma(A_V) = \sigma(A)\setminus\{\lambda\}$, where $A_V$ is the part of the operator $A$ in $V$. Assume that the space $X$ is continuously embedded in a Hilbert space $H$, equipped with a norm $\|\cdot\|_H$ and a scalar product $\<\,\cdot\, , \,\cdot\,\>_H$, such that the operator $A$ has self-adjoint
extension $\widehat{A}:H \supset D(\widehat{A}) \to H$. We say that condition $(G1)$ is satisfied provided
\begin{equation*}
\ \left\{\begin{aligned}
& \text{for any balls } B_1\subset V\cap X^\alpha \text{ and } B_2\subset X_0 \text{ there is } R > 0 \text{ such that } \\
& \<F(t, x + y), x\>_H > - \<F(t, x + y), z\>_H \\
& \text{for } (t,y,z)\in [0,T]\times B_1 \times B_2, \ \text{and} \ x\in X_0 \text{ with } \|x\|_H\ge R.
\end{aligned}\right.
\end{equation*}
Furthermore, we say that condition $(G2)$ is satisfied provided
\begin{equation*}
\ \left\{\begin{aligned}
& \text{for any balls } B_1\subset V\cap X^\alpha \text{ and } B_2\subset X_0 \text{ there is } R > 0 \text{ such that } \\
& \<F(t, x + y), x\>_H < - \<F(t,x + y), z\>_H \\
& \text{for } (t,y,z)\in [0,T]\times B_1 \times B_2, \ \text{and} \ x\in  X_0 \text{ with } \|x\|_H\ge R.
\end{aligned}\right.
\end{equation*}
An important property of these conditions is the fact that if $F$ is a Nemitskii operator associated with $f:[0,+\infty)\times\Omega\times\mathbb{R}\to\mathbb{R}$, then $(G1)$ and $(G2)$ are implicated by well-known Landesman-Lazer conditions introduced in \cite{MR0267269} as well as strong resonance conditions considered in \cite{MR713209}.

Let us observe that the advantage of using the topological degree methods is that, we obtain not only the existence of $T$-periodic solutions, but also we compute its degree. This brings us additional topological information which is useful in the study of the multiplicity and stability of periodic solutions.

The paper is organized as follows.
In Section 2, we provide some spectral properties of the operator ${\bf A}$. In particular we prove that the elements of spectrum with negative real part are actually eigenvalues. The crucial point is Theorem \ref{th-spec-dec}, which express a relationship between spectral decomposition of the operators ${\bf A}$ and $A$. Here the main difficulties are caused by the fact that ${\bf A}$ does not have compact resolvents, despite $A$ has compact resolvents as we assumed.

Section 3 is devoted to the mild solutions for the equation \eqref{row-par} where the nonlinearity depends additionally from a parameter. We provide the standard facts concerning the existence and uniqueness of mild solutions and then, we focus on continuity and compactness properties of the Poincar\'e operator.

In Section 4 we provide geometrical assumptions on the nonlinearity $F$ and use them to prove {\em the degree formula for periodic solutions}, Theorem \ref{th-reso-m-h}, that is the main result of this paper.
Finally, in Section 5 we provide applications of the obtained abstract results to partial differential equations. First of all, in Theorems \ref{lem-est2} and \ref{lem-est3}, we prove that if $F$ is a Nemitskii operator associated with a map $f$, then the well known Landesman-Lazer (see \cite{MR0487001}, \cite{MR0267269}) and strong resonance conditions (see \cite{MR713209}) are actually particular case of $(G1)$ and $(G2)$. As applications we provide criteria on the existence of $T$-periodic solutions for the strongly damped wave equation in terms of Landesman-Lazer and strong resonance conditions. \\[5pt]
\noindent {\em \bf Notation and terminology.} Let $A:X\supset D(A)\to X$ be a linear operator on a real Banach space $X$ equipped with the norm $\|\cdot\|$. The spectrum $\sigma(A)$ of the operator $A$ we define by the complexification method. To be more precise we introduce a new complex vector space $X_\mathbb{C}:=X\times X$ where the operations of addition and multiplication by complex scalar are given by
\begin{equation*}
\begin{aligned}
(x_1,y_1) + (x_2,y_2) & := (x_1 + x_2, y_1 + y_2) && \text{for \ } (x_1, y_1), \ (x_2,y_2) \in X_\mathbb{C}, \\
(\lambda_1 + i\lambda_2)\cdot (x,y) & := (\lambda_1 x - \lambda_2 y, \lambda_1 y + \lambda_2 x) && \text{for \ } (x, y)\in X_\mathbb{C}, \ \lambda_1,\lambda_2\in\mathbb{R}.
\end{aligned}
\end{equation*}
It is well-known that the function
\begin{equation*}
\|z\|_\mathbb{C} := \sup_{\theta\in [0,2\pi]} \|(\sin\theta) x + (\cos\theta) y\| \qquad \mathrm{for} \quad z = (x,y)\in X_\mathbb{C}
\end{equation*}
defines a norm on $X_\mathbb{C}$, and furthermore, $X_\mathbb{C}$ equipped with $\|\cdot\|_\mathbb{C}$ is a Banach space.
From now on we will intuitively denote $x + i y := (x,y)$ for $(x,y)\in X_\mathbb{C}$.

Let us now define the complexification of the linear operator $A$ as a $\mathbb{C}$-linear operator $A_\mathbb{C}$ on $X_\mathbb{C}$ given by the following formula
\begin{align*}
D(A_\mathbb{C}) = D(A)\times D(A), \quad A_\mathbb{C} (x + i y) := A x + i A y \quad\text{for}\quad  x + i y \in D(A_\mathbb{C}).
\end{align*}
By the spectrum $\sigma(A)$ of the operator $A$ we mean the spectrum of its complexification $\sigma(A_\mathbb{C})$. Similarly by the set of eigenvalues $\sigma_p(A)$ of the operator $A$ we understand the set eigenvalues of the operator $A_\mathbb{C}$. For more details on complexification of linear operators we refer the reader to \cite{MR1345385} and \cite{MR1204883}.

We say that the operator $A$ is sectorial provided there are $\phi\in(0,\pi/2)$,  $M\ge 1$ and $a\in\mathbb{R}$, such that the sector $$S_{a,\phi}:=\{\lambda\in\mathbb{C} \ | \ \phi \le |\mathrm{arg} \, (\lambda - a) \le \pi, \ \lambda\neq a\}$$ is contained in the resolvent set of $A$ and furthermore $$\|(\lambda I - A)^{-1}\| \le M/ |\lambda - a| \qquad \mathrm{for} \quad \lambda\in S_{a,\phi}.$$
It is well-known that $-A$ is an infinitesimal generator of analytic semigroup which, throughout this paper, will be denoted by $\{S_A(t)\}_{t\ge 0}$. The operator $A$ is called positive if $\Re \mu > 0$ for any $\mu\in\sigma(A)$. It can be proved that, if $A$ is positive and sectorial, then given $\alpha \ge 0$ the integral
\begin{equation*}
A^{-\alpha} := \frac{1}{\Gamma(\alpha)}\int_0^\infty t^{\alpha - 1}S_A(t) \, dt.
\end{equation*}
is convergent in the uniform operator topology of the space $\mathcal{L}(X)$. Consequently we can define {\em the fractional space} associated with $A$ as the domain of the inverse operator $X^\alpha:= D(A^\alpha)$. The space $X^\alpha$ endowed with the graph norm $\|x\|_\alpha := \|A^\alpha x\|$ is a Banach space, continuously embedded in $X$. For more details on sectorial operators and fractional spaces, we refer the reader to \cite{MR610244}, \cite{MR0089373}, \cite{Pazy}.

\section{Spectral decomposition}

In this section we will assume that $A:X\supset D(A)\to X$ is a sectorial operator on a real Banach space $X$, equipped with the norm $\|\cdot\|$, such that the following conditions are satisfied: \\[3pt]
\noindent\makebox[9mm][l]{$(A1)$}\parbox[t][][t]{118mm}{the operator $A$ is positive and has compact resolvents,}\\[3pt]
\noindent\makebox[9mm][l]{$(A2)$}\parbox[t][][t]{118mm}{there is a Hilbert space $H$ endowed with a scalar product $\<\,\cdot\,, \,\cdot\,\>_H$ and a norm $\|\cdot\|_H$ such that $X$ is embedded in $H$ by a continuous injective map $i:X \hookrightarrow H$,}\\[3pt]
\noindent\makebox[9mm][l]{$(A3)$}\parbox[t][][t]{118mm}{there is a self-adjoint operator $\widehat A:H\supset D(\widehat A) \to H$ such that $\mathrm{Gr}\,(A)\subset \mathrm{Gr}\,(\widehat A)$, where the inclusion is understood in the sense of the product map $i\times i$.}

\begin{remark}\label{rem-pom}
Under the assumptions $(A1)-(A3)$, the spectrum $\sigma(A)$ consists of a sequence (possibly finite) of real positive eigenvalues $$0 < \lambda_1 < \lambda_2 < \ldots < \lambda_i < \lambda_{i+1} < \ldots$$ such that $\dim (\lambda_i I - A) < +\infty$ for $i\ge 1$. \\[2pt] \noindent Indeed, the operator $A$ has compact resolvents which implies that there is a complex sequence $(\lambda_i)_{i\ge 1}$ such that
\begin{equation}\label{nini}
\sigma(A) = \sigma_p(A) = \{\lambda_i \ | \ i\ge 1\} \text{ and } \dim_\mathbb{C}\mathrm{Ker}\,(\lambda_i I - A_\mathbb{C}) < +\infty \quad\text{for} \quad i\ge 1.
\end{equation}
Furthermore the sequence is finite or $|\lambda_i|\to +\infty$ as $n\to +\infty$. Since $\lambda\in\sigma_p(A)$ we have also that $\lambda\in\sigma_p(\widehat A)$ as a consequence of $(A3)$. But $\widehat A$ is a symmetric operator which follows that $\lambda$ is a real number. Since the operator $A$ is positive, the spectrum of $A$ can be exhibited as the increasing sequence of real positive eigenvalues $(\lambda_i)_{i\ge 1}$ which is finite or $\lambda_i \to +\infty$ as $i\to +\infty$. In view of the right part of \eqref{nini} and the fact that $\lambda_i$ are real numbers, we obtain $\dim (\lambda_i I - A) < +\infty$ for $i\ge 1$. \hfill $\square$
\end{remark}

Let us introduce the space ${\bf E}:= X^\alpha\times X$ equipped with the norm
\begin{equation}\label{norm12}
\|(x,y)\|_{\bf E}:= \|x\|_\alpha + \|y\| \qquad \mathrm{for} \quad (x,y)\in {\bf E}.
\end{equation}
We proceed to study the spectral properties of the operator ${\bf A}:{\bf E}\supset D({\bf A})\to{\bf E}$ given by the formula
\begin{equation}
\begin{aligned}\label{op-hyp}
D({\bf A}) & :=\{(x,y)\in {\bf E} = X^\alpha\times X \ | \ x + c y\in D(A)\}, \\
{\bf A}(x,y) & :=(-y, A(x + c y) - \lambda x) \qquad \mathrm{for} \quad (x,y)\in D({\bf A}),
\end{aligned}
\end{equation}
where $\lambda\in\mathbb{R}$ and $c > 0$. If we consider the complexification ${\bf A}_\mathbb{C}:(X^\alpha\times X)_\mathbb{C} \supset D({\bf A}_\mathbb{C}) \to (X^\alpha\times X)_\mathbb{C}$ of the operator ${\bf A}$, then it can be easily checked that ${\bf A}_\mathbb{C}$ is linearly conjugate with the operator ${\bf B}: X_\mathbb{C} \times X^\alpha_\mathbb{C} \supset D({\bf B}) \to X_\mathbb{C} \times X^\alpha_\mathbb{C}$ given by
\begin{equation}
\begin{aligned}\label{nini2}
& D({\bf B}) := \{(x,y)\in X^\alpha_\mathbb{C}\times X_\mathbb{C} \ | \ x + c y \in D(A_\mathbb{C})\}, \\
& {\bf B}(x,y) := (-y, A_\mathbb{C}(x + c y) - \lambda x) \qquad \mathrm{for} \quad (x,y)\in D({\bf B}).
\end{aligned}
\end{equation}
To see this, it is enough to see that $U {\bf A}_\mathbb{C} = {\bf B} U$, where $U:(X^\alpha\times X)_\mathbb{C} \to X^\alpha_\mathbb{C}\times X_\mathbb{C}$ is a $\mathbb{C}$-linear isomorphism given by
$$U((x_1, y_1) + i(x_2, y_2)) := (x_1 + i x_2, y_1 + i y_2) \ \text{ for } \ (x_1, y_1) + i(x_2, y_2)\in (X^\alpha\times X)_\mathbb{C}.$$ Hence, without loss of generality we will implicitly consider the operator ${\bf B}$ instead of ${\bf A}_\mathbb{C}$. In the following theorem we provide description of the spectrum of ${\bf A}$.
\begin{theorem}\label{th-spec}
The following assertions hold.
\begin{enumerate}
\item[(i)] The set $\sigma({\bf A}) \setminus\{1/c\}$ consists of the eigenvalues of the operator ${\bf A}$. \\[-8pt]
\item[(ii)] If $\lambda_k \le \lambda < \lambda_{k+1}$ for some $k\ge 1$, then $$\{\mu\in\sigma_p({\bf A}) \ | \ \Re\mu\le 0\} = \{\mu_i^- \ | \ 1\le i\le k\},$$ where
\begin{equation}\label{l1}
\mu_i^\pm := \frac{\lambda_ic \pm \sqrt{(\lambda_ic)^2 - 4(\lambda_i - \lambda)}}{2} \ \ \text{ for }\ i\ge 1.
\end{equation}
Furthermore, if $\lambda < \lambda_1$ then $\{\mu\in\sigma_p({\bf A}) \ | \ \Re\mu\le 0\} = \emptyset$. \\[-8pt]
\item[(iii)] If $\lambda_k \le \lambda < \lambda_{k+1}$ for some $k\ge 1$ and $\Re\mu \le 0$, then
$$\mathrm{Ker}\,(\mu  I - {\bf A}) = \mathrm{Ker}\,(e^{-\mu t} I - S_{\bf A}(t)) \qquad \mathrm{for} \quad t > 0.$$
\end{enumerate}
\end{theorem}
In the proof of the theorem we will use the following lemma.
\begin{lemma}\label{lem-11}
Let $\mu\in\mathbb{C}\setminus\{1/c\}$ be arbitrary. Then $(x,y)\in \mathrm{Ker}\,(\mu I - {\bf A}_\mathbb{C})$ if and only if $(x,y) = (w,-\mu w)$ for some $w\in \mathrm{Ker}\,\left((\lambda - \mu^2)/(1 - c\mu) I -A_\mathbb{C}\right)$.
\end{lemma}
\noindent\textbf{Proof.} Assume that $\mu\in\mathbb{C}\setminus\{1/c\}$ and $\mu (x,y) = {\bf A}_\mathbb{C}(x,y)$ for some $(x,y)\neq 0$. In view of \eqref{nini2}, it implies that $x + c y\in D(A_\mathbb{C})$ and
\begin{equation*}
\mu x = - y, \quad \mu y = A_\mathbb{C}(x + c y) - \lambda x.
\end{equation*}
Hence $x\neq 0$, $(1 - c\mu)x\in D(A_\mathbb{C})$ and $(\lambda - \mu^2)x = A_\mathbb{C}((1 - c\mu)x)$, which gives $$(\lambda - \mu^2)/(1 - c\mu)x = A_\mathbb{C} x \quad\text{and} \quad (\lambda - \mu^2)/(1 - c\mu) = \lambda_l \ \text{ for some }l\ge 1,$$ where $\lambda_l$ is the $l$-th eigenvalue of $A$ (see Remark \ref{rem-pom}).
Now if we take $w:=x$, then $x\in\mathrm{Ker}\,(\lambda_l I - A_\mathbb{C})$ and $(x,y) = (x, - \mu x) = (w,-\mu w)$, which completes "if" part of the proof. On the other hand, if $(x,y) = (w, - \mu w)$ for some
$w\in \mathrm{Ker}\,\left((\lambda - \mu^2)/(1 - c\mu) I - A_\mathbb{C}\right)$, then $x+cy\in D(A_\mathbb{C})$ and
\begin{align*}
\mu y - A_\mathbb{C}(x + c y) + \lambda x & = (\lambda - \mu^2)w - A_\mathbb{C}((1-c\mu)w) \\
& = (\lambda - \mu^2)w - (1-c\mu)A_\mathbb{C} w =0.
\end{align*}
Consequently $\mu (x,y) = {\bf A}_\mathbb{C}(x,y)$, which gives desired conclusion. \hfill $\square$\\[5pt]

\noindent\textbf{Proof of Theorem \ref{th-spec}.}  For the point $(i)$, take $\mu\in\sigma({\bf A}) \setminus\{1/c\}$ and suppose that $\mathrm{Ker}\,(\mu I - {\bf A}_\mathbb{C})=\{0\}$. We claim that for every $(f,g)\in X^\alpha_\mathbb{C}\times X_\mathbb{C}$ there is $(x,y)\in D({\bf A}_\mathbb{C})$ such that $\mu (x,y) - {\bf A}_\mathbb{C}(x,y) = (f,g)$. Assume for the moment that our claim is true. Then the inverse operator $(\mu I - {\bf A}_\mathbb{C})^{-1}$ is bounded on $X^\alpha\times X$ because ${\bf A}_\mathbb{C}$ is closed. Therefore $\mu\in\varrho({\bf A})$, which contradicts the fact that $\mu$ belongs to spectrum and proves that $\mathrm{Ker}\,(\mu  I - {\bf A}_\mathbb{C})\neq\{0\}$. Now we proceed to the proof of the claim. Let us take $(f,g)\in X^\alpha_\mathbb{C}\times X_\mathbb{C}$ and consider the following equations
\begin{equation}\label{rro}
\mu x = - y + f, \quad \mu y =  A_\mathbb{C}(x + c y) - \lambda x + g.
\end{equation}
Multiplying the former equation by $c\lambda - \mu$ and later by $1 - \mu c$, one has
\begin{align*}
(c\lambda - \mu)\mu x & = -(c\lambda - \mu)y + (c\lambda - \mu)f \\
(1 - \mu c)\mu y & = (1 - \mu c)A_\mathbb{C}(x + c y) - (1 - \mu c)\lambda x + (1 - \mu c)g.
\end{align*}
Summing up these equations we obtain
\begin{equation*}
    (\lambda - \mu^2) (x + cy) = (1 - \mu c)A_\mathbb{C}(x + c y) + (1 - \mu c)h,
\end{equation*}
where $h:=(c\lambda - \mu)/(1 - c\mu) f + g$, and hence
\begin{equation*}
    \frac{\lambda - \mu^2}{1 - c\mu} (x + cy) = A_\mathbb{C}(x + c y) + h.
\end{equation*}
Since $\mu\neq 1/c$ and $\mathrm{Ker}\,(\mu I - {\bf A}_\mathbb{C})=\{0\}$, from the Lemma \ref{lem-11} it follows that $(\lambda - \mu^2)/(1 - c\mu)$ is an element from the resolvent set of the operator $A_\mathbb{C}$. Therefore
\begin{equation*}
    x + c y = b:= \left(\frac{\lambda - \mu^2}{1 - \mu c} I - A_\mathbb{C}\right)^{-1} h,
\end{equation*}
which allows us to define
$$x := \frac{1}{1 - \mu c}(b - c f) \quad\text{and}\quad y := \frac{1}{1 - \mu c}(f - \mu b).$$
Since $b\in D(A_\mathbb{C})\subset X^\alpha_\mathbb{C}$ and $f\in X^\alpha_\mathbb{C}$, one has $(b - c f)\in X^\alpha_\mathbb{C}$ and hence $x\in X^\alpha_\mathbb{C}$. Furthermore $x + c y = b\in D(A_\mathbb{C})$, which implies that $(x,y)\in D({\bf A}_\mathbb{C})$. Therefore, it is enough to check that the equations \eqref{rro} are satisfied. To this end observe that
$$\mu x = \frac{\mu}{1 - \mu c}(b - c f) = - \frac{1}{1 - \mu c}(f - \mu b) + f = - y + f,$$
and furthermore, we have the following sequence of equivalent equalities
\begin{align*}
& \mu y = A_\mathbb{C}(x + c y) - \lambda x + g \\
& \frac{\mu}{1 - \mu c}(f - \mu b) = A_\mathbb{C} b - \frac{\lambda}{1 - \mu c}(b - c f) + g \\
& \frac{\lambda - \mu^2}{1 - \mu c} b - A_\mathbb{C} b = \frac{c\lambda - \mu}{1 - \mu c} f + g \\
& h = \frac{c\lambda - \mu}{1 - \mu c} f + g.
\end{align*}
The last equality is true by the definition of $h$ and hence, the assertion $(i)$ follows. \\[5pt]
To verify $(ii)$ let us take $\mu\in\sigma_p({\bf A})$ such that $\Re\mu\le 0$. Then, by Lemma \ref{lem-11} and Remark \ref{rem-pom}, one has $(\lambda - \mu^2)/(1 - c\mu) = \lambda_i$ for some $i\ge 1$. Hence the equation $\mu^2 - c\lambda_i\mu + \lambda_i - \lambda = 0$ is satisfied and computing its roots we infer that either $\mu = \mu_i^+$ or $\mu = \mu_i^-$. Since $\lambda_k \le \lambda < \lambda_{k+1}$ and $\Re \mu \le 0$, it follows that $\mu = \mu_i^-$ for some $1\le i\le k$ and hence $\{\mu\in\sigma_p({\bf A}) \ | \ \Re\mu\le 0\} \subset \{\mu_i^- \ | \ 0\le i\le k\}$. In order to prove the opposite inclusion take $\mu = \mu_i^-$ for $1 \le i\le k$. Then $\mu \le 0$ and the equation $(\lambda - (\mu_i^-)^2)/(1 - c\mu_i^-) = \lambda_i$ is satisfied. Therefore Lemma \ref{lem-11} shows that $\mu_i^-\in \sigma_p({\bf A})$, which completes the proof of $(ii)$. \\[5pt]
To see $(iii)$ observe that from the point $(ii)$ it follows that the set $$\{\tilde\mu\in\sigma_p({\bf A}) \ | \ \Re\tilde\mu = \Re\mu \}$$ consists of at most one eigenvalue. To complete the proof, it remains to observe that Theorem 16.7.2 from \cite{MR0089373} leads to
\begin{equation*}
\mathrm{Ker}\,(e^{-\mu t} I - S_{\bf A}(t)) = \overline{\mathrm{Ker}\,(\mu  I - {\bf A})} = \mathrm{Ker}\,(\mu  I - {\bf A}) \qquad \mathrm{for} \quad t>0,
\end{equation*}
where the last equality follows from the fact that the operator $\mu  I - {\bf A}$ is closed and hence its kernel is closed as well. This completes the proof of theorem. \hfill $\square$ \\

Let us now proceed to study the spectral decomposition of the operator ${\bf A}$. For this purpose we will need the following proposition.

\begin{proposition}{\em (\cite[Theorem 2.3]{Kok2})}\label{th:10}
Let us denote by $\lambda := \lambda_k$ the $k$-th eigenvalue of $A$ and let $X_0 := \mathrm{Ker}\, (\lambda I - A)$. Then there are closed subspaces $X_+$, $X_-$ of $X$ such that $X = X_+\oplus X_-\oplus X_0$ and the following assertions hold. \\[-10pt]
\begin{enumerate}
\item[(i)] We have the following inclusions $$X_-\subset D(A), \ \ A(X_-)\subset X_- \text{ \ \ and \ \ } A(X_+\cap D(A)) \subset X_+.$$ Furthermore $X_-$ is a finite dimensional and $X_-= \{0\}$ if $k=1$ and $$X_-=\mathrm{Ker}\,(\lambda_1 I - A)\oplus\ldots\oplus\mathrm{Ker}\,(\lambda_{k-1} I - A) \text{ \ if \ } k\ge 2.$$
\item[(ii)] If $A_+:X_+\supset D(A_+) \to X_+$ and $A_-:X_-\supset D(A_-) \to X_-$ are parts of the operator $A$ in the spaces $X_+$ and $X_-$, respectively, then $$\qquad\qquad\sigma(A_+) = \{\lambda_i \ | \ i\ge k+1 \} \text{ \ and \ } \sigma(A_-) = \{\lambda_i \ | \ 1\le i \le k-1 \} \text{ \ for \ }k\ge 1.$$
\item[(iii)] The spaces $X_0$, $X_-$ are $X_+$ mutually orthogonal in $H$, which means that $\<i (u_l),i(u_m)\>_H = 0$ for $u_l\in X_l$, $u_m\in X_m$ where $l,m\in\{0,-,+\}$, $l\neq m$. \\[-8pt]
\end{enumerate}
\end{proposition}

\begin{remark}\label{rem-proj}
Consider the decomposition $X = X_+\oplus X_-\oplus X_0$ obtained in Proposition \ref{th:10}. Let us denote by $P, Q_\pm:X\to X$ projections given by
\begin{equation}\label{wzz1ab}
    P x = x_0 \ \text{ and } \ Q_\pm x = x_\pm \qquad \mathrm{for} \quad x\in X,
\end{equation}
where $x = x_+ + x_0 + x_-$ for $x_i\in X_i$, $i\in \{0,-,+\}$. Since the components are closed in $X$, the projections are continuous maps. Let us denote $Q:= Q_- + Q_+$. Since the inclusion $X^\alpha\subset X$ is continuous, one can decompose $X^\alpha$ on a direct sum of closed spaces
$$X^\alpha = X_0\oplus X^\alpha_-\oplus X^\alpha_+, \text{ \ where \ } X^\alpha_-:= X^\alpha\cap X_-, \ X^\alpha_+:=X^\alpha\cap X_+.$$
Therefore the projections $P$ and $Q_\pm$ can be also considered as continuous maps $P, Q_\pm:X^\alpha\to X^\alpha$, given for any $x\in X^\alpha$ by the formula \eqref{wzz1ab}.
\end{remark}
Let us now proceed to the main result of this section. The following theorem express a relationship between spectral decomposition of the operators ${\bf A}$ and $A$.
\begin{theorem}\label{th-spec-dec}
Let us denote by $\lambda := \lambda_k$ the $k$-th eigenvalue of $A$ and let ${\bf E}_0 := \mathrm{Ker}\,(\lambda I - A)\times\mathrm{Ker}\,(\lambda I - A)$. Then there are closed subspaces ${\bf E}_+,{\bf E}_-$ of ${\bf E}$ such that ${\bf E}:={\bf E}_-\oplus{\bf E}_0\oplus{\bf E}_+$ and the following assertions hold.
\begin{enumerate}
\item[(i)] We have the following inclusions $${\bf E}_-\subset D({\bf A}), \ \ {\bf A}({\bf E}_-)\subset {\bf E}_- \text{ \ \ and \ \ } {\bf A}({\bf E}_+\cap D({\bf A})) \subset {\bf E}_+.$$
Furthermore, if $k\ge 2$ then ${\bf E}_- = K_1^-\oplus K_2^-\oplus\ldots\oplus K_{k-1}^-$, where $$K_i^- := \{(w,-\mu_i^- w) \in{\bf E}\ | \ w\in\mathrm{Ker}\,(\lambda_i I - A)\}$$
and $\mu_i^-$ for $i\ge 1$ are numbers given by \eqref{l1}. If $k=1$ then ${\bf E}_- = \{0\}$.  \\[-7pt]
\item[(ii)] If ${\bf A}_+:{\bf E}_+\supset D({\bf A}_+)\to{\bf E}_+$ and ${\bf A}_-:{\bf E}_-\supset D({\bf A}_-)\to{\bf E}_-$ are parts of ${\bf A}$ in the spaces ${\bf E}_-$ and ${\bf E}_+$, respectively, then $$\qquad\quad\sigma({\bf A}_+)\subset \{z\in\mathbb{C} \ | \ \Re z > 0\} \text{ \ and \ }\sigma({\bf A}_-) = \{\mu_i^- \ | \ 1\le i \le  k - 1\}.$$
\item[(iii)] If ${\bf P},{\bf Q}_-,{\bf Q}_+:{\bf E}\to{\bf E}$ are projections on the spaces ${\bf E}_0$, ${\bf E}_-$ and ${\bf E}_+$, respectively, and ${\bf Q} := {\bf Q}_- + {\bf Q}_-$, then
\begin{align}\label{proj}
{\bf P}(x,y) = (Px, Py) \ \text{ and } \ {\bf Q}(x,y) = (Qx, Qy) \ \text{ for } \ (x,y)\in{\bf E}.
\end{align}
\end{enumerate}
\end{theorem}
In the proof we will use the following two lemmata.
\begin{lemma}\label{lem-22}
Assume that $\mu_i^\pm$ for $i\ge 1$ are numbers given by \eqref{l1}. If $1\le i \le k$ then $\mu_i^\pm$ is real and the subspaces
\begin{equation}\label{nini3}
K_i^\pm := \{(w,-\mu_i^\pm w) \in{\bf E}\ | \ w\in\mathrm{Ker}\,(\lambda_i I - A)\}
\end{equation}
are such that $K_i^\pm \subset \mathrm{Ker}\,(\mu_i^\pm  I - {\bf A})$ and $\mathrm{Ker}\, (\lambda_i I - A)\times \mathrm{Ker}\, (\lambda_i I - A) = K_i^+\oplus K_i^-$.
\end{lemma}
\noindent\textbf{Proof.} Let us first observe that $\mu_i^+ \neq \mu_i^-$ for $1\le i \le k$, which implies that $K_i^+\cap K_i^- = \{0\}$. By Remark \ref{rem-pom}, one has $\dim\mathrm{Ker}\,(\lambda_i I - A) < + \infty$ and $$\dim \mathrm{Ker}\, (\lambda_i I - A)\times \mathrm{Ker}\, (\lambda_i I - A) = 2 \dim\mathrm{Ker}\,(\lambda_i I - A) = \dim K_i^+ + \dim K_i^-.$$ Hence it follows that $\mathrm{Ker}\, (\lambda_i I - A)\times \mathrm{Ker}\, (\lambda_i I - A) = K_i^+\oplus K_i^-$. To complete the proof of theorem let us take $(x,y)\in K_i^\pm$. Then $(x,y) = (w,-\mu_i^\pm w)$ for some $w\in\mathrm{Ker}\,(\lambda_i I - A)$, which implies that $x + cy \in D(A)$ and
\begin{align*}
\mu_i^\pm y - A(x + c y) + \lambda x & = (\lambda - (\mu_i^\pm)^2)w - A((1-c\mu_i^\pm)w) \\
& = (\lambda - (\mu_i^\pm)^2)w - (1-c\mu_i^\pm)A w \\
& = (\lambda - (\mu_i^\pm)^2)w - (1-c\mu_i^\pm)\lambda_i w \\
& = - ((\mu_i^\pm)^2 - \lambda_ic \mu_i^\pm + \lambda_i - \lambda) w = 0.
\end{align*}
The last equality follows from the fact that $\mu_i^\pm$ are the roots of the equation $\mu^2 - \lambda_ic \mu + \lambda_i - \lambda = 0$. Hence $\mu_i^\pm (x,y) = (- y, A(x + c y) - \lambda x) = {\bf A} (x,y)$, which implies that
$K_i^\pm\subset \mathrm{Ker}\,(\mu_i^\pm  I - {\bf A})$ as desired. \hfill $\square$ \\

\begin{lemma}\label{lem-wart}
Assume that $B:V\to V$ is a linear operator on a real finite dimensional space $V$ and let $V_1,V_2,\ldots,V_l$ be subspaces such that $V=V_1\oplus V_2\oplus\ldots\oplus V_l$. If the real numbers $(\nu_i)_{i=1}^l$ are such that $B x = \nu_i x$ for $x\in V_i$, then
\begin{enumerate}
\item[(a)] $\sigma(B) = \{\nu_i \ | \ 1\le i \le l\}$, \\[-10pt]
\item[(b)] for any $1\le i \le l$ one has $N_{\nu_i}(B):=\bigcup_{m=1}^\infty \mathrm{Ker}\,(\nu_i I - B)^m = \mathrm{Ker}\,(\nu_i I - B)$.
\end{enumerate}
\end{lemma}
\noindent\textbf{Proof.} $(a)$ It is enough to prove that $\sigma(B) \subset \{\nu_i \ | \ 1\le i \le l\}$. The opposite inclusion is obvious. Let $\nu\in\mathbb{C}$ be such that
$\nu z = B_\mathbb{C} z$ for some $z := x + iy\in V_\mathbb{C} $, $z\neq 0$. Then $V_\mathbb{C} = (V_1\times V_1)\oplus (V_2\times V_2)\oplus\ldots\oplus (V_l\times V_l)$ and $z=z_1+z_2+\ldots+ z_l$ where $z_i\in V_i\times V_i$ for $1\le i\le l$. Therefore
$\nu z = B_\mathbb{C} z =  \nu_1 z_1 + \nu_2 z_2 + \ldots + \nu_l z_l$. Since $z\neq 0$, there is $1\le i \le l$ such that $z_i\neq 0$ and consequently $\nu = \nu_i$, which proves desired inclusion. \\[5pt]
$(b)$ It suffices to prove that $N_{\nu_i}(B) \subset \mathrm{Ker}\,(\nu_i I - B)$. Take $z\in N_{\nu_i}(B)\setminus \{0\}$. Then there is $i_0 \ge 1$ such that $(\nu_i I - B)^{i_0} z = 0$ and there are $z_i\in V_i$ ($1\le i\le l$) such that $z = z_1 + z_2 + \ldots + z_l$. Hence
\begin{align*}
0 & = (\nu_i I - B)^{i_0} z = (\nu_i I - B)^{i_0} z_1 + (\nu_i I - B)^{i_0} z_2 + \ldots + (\nu_i I - B)^{i_0} z_l \\
& = (\nu_i - \nu_1)^{i_0} z_1 + (\nu_i - \nu_2)^{i_0} z_2 + \ldots + (\nu_i - \nu_l)^{i_0} z_l.
\end{align*}
Since $z \neq 0$, one of the elements $z_1, z_2, \ldots, z_n$ is also nonzero. Suppose that $z_j \neq 0$ for some $1\le j \le l$. Then $(\nu_i - \nu_j)^{i_0} z_j = 0$ and hence $\nu_i = \nu_j$. It implies that $z\in \mathrm{Ker}\,(\nu_i I - B)$, which completes the proof of desired inclusion. \hfill $\square$ \\

\noindent\textbf{Proof of Theorem \ref{th-spec-dec}.} If $k = 1$ then let us define
\begin{align*}
{\bf E}_- := \{0\}, \quad {\bf E}_+ := (X^\alpha\cap X_+)\times X_+.
\end{align*}
For $k\ge 2$ we make use of \eqref{nini3} and write
$$M_1:= K_1^+\oplus\ldots\oplus K_{k-1}^+ \quad\text{and}\quad M_2:=(X^\alpha\cap X_+)\times X_+.$$ Then let us define
\begin{align*}
{\bf E}_- := K_1^-\oplus \ldots\oplus K_{k-1}^-, \quad {\bf E}_+ := M_1\oplus M_2.
\end{align*}
It is not difficult to check that ${\bf E}_0$, ${\bf E}_-$, ${\bf E}_+$ are closed subspaces such that ${\bf E}={\bf E}_-\oplus{\bf E}_0\oplus{\bf E}_+$ and the assertion $(iii)$ holds.

Observe that ${\bf E}_0\subset D({\bf A})$ and ${\bf A}({\bf E}_0)\subset {\bf E}_0$. Furthermore, from Lemma \ref{lem-22} it follows that ${\bf E}_-\subset D({\bf A})$ and ${\bf A}({\bf E}_-)\subset {\bf E}_-$, because $K_i^-$ is contained in the eigenspace of the operator ${\bf A}$ for $1\le i\le k-1$. For the proof of $(i)$, it remains to verify the inclusion ${\bf A}(D({\bf A})\cap {\bf E}_+)\subset {\bf E}_+$. Similarly as before Lemma \ref{lem-22} says that, for any $1\le i \le k-1$, the elements of $K_i^+$ are contained in the eigenspace of ${\bf A}$ and therefore ${\bf A}(M_1) \subset M_1$. It remain to show that ${\bf A}(D({\bf A})\cap M_2) \subset M_2$. To this end take $(x,y)\in D({\bf A})\cap M_2$. Then $x + c y\in D(A)$ and $x\in X^\alpha$, so that $y\in X^\alpha$, which yields $y\in X^\alpha\cap X_+$. Observe that $A(x + c y) - \lambda x\in X_+$, which is a consequence of the fact that $x + cy \in D(A)\cap X_+$ and the inclusion $A(D(A)\cap X_+)\subset X_+$ from Proposition \ref{th:10} $(i)$. Therefore ${\bf A}(x,y)\in M_2$ as desired and proof of $(i)$ is completed. \\[5pt]
For the proof of point $(ii)$ let us first show that $\Re\mu > 0$ for $\mu\in\sigma({\bf A}_+)$.

\noindent Suppose to the contrary that $\mu\in\sigma({\bf A}_+)$ and $\Re\mu \le 0$. Let ${\bf A}_+^1$ and ${\bf A}_+^2$ be the parts of the operator ${\bf A}_+$ in $M_1$ and $M_2$, respectively. By inclusion ${\bf A}(M_1) \subset M_1$ and Lemma \ref{lem-wart}, it follows that the spectrum of $\mu  I - {\bf A}_+^1$ consists of its eigenvalues $\{\mu - \mu_i^+ \ | \ 1\le i \le k-1\}$. Since $\Re\mu \le 0$ and $\mu_i^+ > 0$ for $1\le i \le k$, we deduce that the complexification of the operator $\mu  I - {\bf A}_+^1$ is bijection and hence $\mu\in\sigma({\bf A}_+^2)$. Observe that the operator ${\bf A}_+^2$ can be given by the formula
\begin{align*}
& D({\bf A}_+^2) = \{(x,y)\in (X^\alpha\cap X_+)\times X_+ \ | \ x + c y \in D(A_+)\}, \\
& {\bf A}_+^2(x,y) = (-y, A_+(x + c y) - \lambda x) \qquad \mathrm{for} \quad (x,y)\in D({\bf A}_+^2).
\end{align*}
By Proposition \ref{th:10} $(ii)$ one has $\sigma(A_+) = \{\lambda_i \ | \ i\ge k+1 \}$. Hence we can apply Theorem \ref{th-spec} to the operator $A_+$ and derive that the set $\sigma({\bf A}_+^2) \setminus\{1/c\}$ consists of eigenvalues and $\{\tilde\mu\in\sigma_p({\bf A}_+^2) \ | \ \Re\tilde\mu\le 0\}=\emptyset$. This is a contradiction because $\Re\mu\le 0$ and $\mu\in\sigma({\bf A}_+^2)$. This yields the first inclusion of $(ii)$. On the other hand, the equality $\sigma({\bf A}_-) = \{\mu_1^-, \mu_2^-, \ldots , \mu_{k-1}^-\}$ is a consequence of the inclusion ${\bf A}({\bf E}_-)\subset {\bf E}_-$ and Lemmata \ref{lem-22} and \ref{lem-wart}. Hence the proof of $(ii)$ is completed. \hfill $\square$ \\

Let us now remind that ${\bf A}$ is a sectorial operator (see for example \cite{MR2571574}, \cite{MR702424}) and hence $-{\bf A}$ generates a $C_0$ semigroup $\{S_{\bf A}(t)\}_{t\ge 0}$ of bounded operators on ${\bf E}$. Properties of this semigroup, that will be used in this paper, are contained in the following corollary that is a consequence of Theorem \ref{th-spec-dec} and \cite[Theorem 1.5.3]{MR610244}.

\begin{corollary}\label{cor-dec-2aba}
Let us denote by $\lambda := \lambda_k$ the $k$-th eigenvalue of the operator $A$ and let ${\bf E}={\bf E}_-\oplus{\bf E}_0\oplus{\bf E}_+$ be the direct sum decomposition obtained in Theorem \ref{th-spec-dec}.
\begin{itemize}
\item[(i)] We have the following inclusions $S_{\bf A}(t) {\bf E}_i \subset {\bf E}_i$ for $t\ge 0$, $i\in\{0,-,+\}$.
In particular, for any $t\ge 0$ and $z\in {\bf E}$
\begin{equation}\label{ds1-hyp}
S_{\bf A}(t){\bf P} z = {\bf P} S_{\bf A}(t)z \quad\text{and}\quad S_{\bf A}(t){\bf Q}_\pm z = {\bf Q}_\pm S_{\bf A}(t) z.
\end{equation}
\item[(ii)] The $C_0$ semigroup $\{S_{\bf A}(t)|_{{\bf E}_-}\}_{t\ge 0}$ can be uniquely extended to a $C_0$ group on ${\bf E}_-$ and there are constants $M, \delta > 0$ such that
\begin{alignat}{2}\label{nieee1}
\|S_{\bf A}(t)z\|_{\bf E} & \le M e^{- \delta t} \|z\|_{\bf E} && \quad\text{for}\quad z\in {\bf E}_+, \ t\ge 0, \\ \label{nieee2}
\|S_{\bf A}(t)z\|_{\bf E} & \le M e^{ \delta t} \|z\|_{\bf E} && \quad\text{for}\quad z\in {\bf E}_-, \ t\le 0.
\end{alignat}
\item[(iii)] If ${\bf A}_0:{\bf E}_0 \to {\bf E}_0$ is a part of ${\bf A}$ in ${\bf E}_0$, then
\begin{equation}
\begin{aligned}\label{dss1}
{\bf A}_0(x,y) = (-y, c\lambda y) \text{ \ \ and \ \ } S_{{\bf A}_0}(t)(x,y) = S_{\bf A}(t)(x,y)
\end{aligned}
\end{equation}
for $(x,y)\in{\bf E}_0$ and $t\ge 0$.
\end{itemize}
\end{corollary}

\section{Continuity and compactness properties of Poincar\'e operator}

Assume that $A: X\supset D(A)\to X$ is a positively defined sectorial operator with compact resolvents on a separable Banach space $X$, equipped with a norm $\|\cdot\|$. Consider the following family of differential equations
\begin{equation}\label{rownn}
    \ddot u(t) = -A u(t) - c A \dot u(t) + \lambda u(t) + F(s, t, u(t)), \qquad t\ge 0
\end{equation}
where $\lambda$ is a real number, $c > 0$, $s\in[0,1]$ is a parameter and $F:[0,1]\times [0,+\infty)\times X^\alpha\to X$ is a continuous map satisfying the following conditions: \\[2pt]
\noindent\makebox[9mm][l]{$(F1)$}\parbox[t]{118mm}{for every $x\in X^\alpha$ there is an open neighborhood $V\subset X^\alpha$ of $x$ and constant $L > 0$ such that for any $s\in[0,1]$, $x_1,x_2\in V$ and $t\in [0,+\infty)$ $$\|F(s,t,x_1) - F(s,t,x_2)\|\le L \|x_1 - x_2\|_\alpha;$$}\\
\noindent\makebox[9mm][l]{$(F2)$}\parbox[t]{118mm}{there is a continuous function $c:[0,+\infty)  \to [0,+\infty)$ such that $$\|F(s,t,x)\| \le c(t)(1 + \|x\|_\alpha)\quad \mbox{ for } \  s\in[0,1], \ t\in [0,+\infty), \ x\in X^\alpha;$$}\\
\noindent\makebox[9mm][l]{$(F3)$}\parbox[t][][t]{118mm}{if $V\subset X^\alpha$ is a bounded set, then $F([0,1]\times [0,+\infty)\times V)$ is a relatively compact subset of $X$.}\\[5pt]
Let us note that the equation \eqref{rownn} can be written in the form
\begin{equation}\label{row-abstr}
    \dot w(t) = -{\bf A} w(t) + {\bf F}(s,t, w(t)), \qquad t > 0
\end{equation}
where ${\bf A}:{\bf E}\supset D({\bf A})\to{\bf E}$ is a linear operator given by \eqref{op-hyp} and ${\bf F}:[0,1]\times [0,+\infty)\times {\bf E}\to{\bf E}$ is a map defined by
\begin{equation*}
    {\bf F}(s, t, (x,y)):=(0,F(s, t, x)) \ \ \text{for} \ s\in [0,1], \ t\in[0,+\infty), \ (x,y)\in{\bf E}.
\end{equation*}
\begin{definition}
We say that a continuous map $w:[0, +\infty) \to {\bf E}$ is \emph{a (global) mild solution} of equation \eqref{row-abstr} starting at $z_0\in{\bf E}$, provided
\begin{equation*}
w(t) = S_{\bf A}(t)z_0 + \int_{0}^t S_{\bf A}(t - \tau){\bf F}(s,\tau,w(\tau)) \, d\tau \qquad \mathrm{for} \quad t\ge 0.
\end{equation*}
\end{definition}
Let us now proceed to the following proposition, which will be helpful in the study of compactness properties of mild solutions for the equation \eqref{row-abstr}.
\begin{proposition}\label{lem-con-sem}
There is a norm $|\cdot|$ on ${\bf E}$, equivalent with $\|\cdot\|_{\bf E}$, such that:
\begin{enumerate}
\item[(a)] the following inequality holds
\begin{equation}\label{rnon1ab}
     |{\bf Q}_+z| \le |z| \qquad \mathrm{for} \quad z\in {\bf E},
\end{equation}
where ${\bf Q}_+$ is the projection obtained in Theorem \ref{th-spec-dec};\\[-8pt]
\item[(b)] for any bounded $\Omega\subset{\bf E}$ one has
\begin{equation}\label{rownorm}
\beta(S_{\bf A}(t)\Omega) \le e^{-\delta t}\beta(\Omega) \qquad \mathrm{for} \quad t\ge 0,
\end{equation}
where $\beta$ is the Hausdorff measure of noncompactness associated with $|\cdot|$ \footnote{\, See Appendix \ref{mesu}} and $\delta > 0$ is a constant from the inequality \eqref{nieee1}.
\end{enumerate}
\end{proposition}
\noindent\textbf{Proof.} Let us take the remaining two projection ${\bf P}$ and ${\bf Q}_-$ obtained in Theorem \ref{th-spec-dec} and denote ${\bf P}_- := {\bf P} + {\bf Q}_- $. Using the Inverse Mapping Theorem it is not difficult to see that $\|\cdot\|_{\bf E}$ is equivalent with the norm given by $$\|z\|_1:= \|{\bf P}_- z\|_{\bf E} + \|{\bf Q}_+ z\|_{\bf E} \qquad \mathrm{for} \quad z\in {\bf E},$$ and hence there are constants $c_1,c_2 > 0$ such that $c_1\|z\|_{\bf E} \le \|z\|_1 \le c_2\|z\|_{\bf E}$ for $z\in {\bf E}$. Let us define the norm
$$|z| := \|{\bf P}_- z\|_{\bf E} + \sup_{t\ge 0}\|e^{\delta t}S_{\bf A}(t){\bf Q}_+ z\|_{\bf E} \qquad \mathrm{for} \quad z\in {\bf E},$$
where $\delta > 0$ is a constant from the inequality \eqref{nieee1}.
It is not difficult to check that $$c_1\|z\|_{\bf E} \le \|z\|_1 \le |z| \qquad \mathrm{for} \quad z\in {\bf E}$$ and furthermore, by \eqref{nieee1}, one has $$|z| \le \|{\bf P}_- z\|_{\bf E} + M \|{\bf Q}_+ z\|_{\bf E}  \le M \|z\|_1 \le c_2 M \|z\|_{\bf E} \qquad \mathrm{for} \quad z\in {\bf E}.$$ Hence the norms $|\cdot|$ and $\|\cdot\|_{\bf E}$ are equivalent. For the proof of \eqref{rnon1ab} observe that
\begin{align*}
|{\bf Q}_+z| & = \|{\bf P}_- {\bf Q}_+z\|_{\bf E} + \sup_{t\ge 0}\|e^{\delta t}S(t){\bf Q}_+^2 z\|_{\bf E} \\
& = \sup_{t\ge 0}\|e^{\delta t}S(t){\bf Q}_+ z\|_{\bf E} \le |z|  \quad \text{for} \quad z\in {\bf E}
\end{align*}
and point $(a)$ follows. To prove point $(b)$ let us take $t\ge 0$ and $z\in {\bf E}_+$. Then from \eqref{ds1-hyp} it follows that
\begin{align*}
|S_{\bf A}(t)z| & = \|{\bf P}_- S_{\bf A}(t) z\|_{\bf E} + \sup_{s\ge 0}\|e^{\delta s}S_{\bf A}(s){\bf Q}_+ S_{\bf A}(t)z\|_{\bf E} \\
& = \|S_{\bf A}(t) {\bf P}_- z\|_{\bf E} + \sup_{s\ge 0}\|e^{\delta s}S_{\bf A}(s) S_{\bf A}(t){\bf Q}_+ z\|_{\bf E} \\
& = e^{-\delta t} \sup_{s\ge 0}\|e^{\delta (t + s)}S_{\bf A}(t + s){\bf Q}_+ z\|_{\bf E} \le e^{-\delta t}  \sup_{s\ge 0}\|e^{\delta s}S_{\bf A}(s){\bf Q}_+ z\|_{\bf E} \\
& \le e^{-\delta t} \left(\|{\bf P}_- z\|_{\bf E} + \sup_{s\ge 0}\|e^{\delta s}S_{\bf A}(s){\bf Q}_+ z\|_{\bf E}\right) = e^{-\delta t} |z|.
\end{align*}
Consequently, we infer that
\begin{equation}\label{rnon1}
    |S_{\bf A}(t)z| \le e^{-\delta t} |z| \qquad \mathrm{for} \quad t\ge 0 \text{ \ and \ } z\in {\bf E}_+.
\end{equation}
By the properties of the measure $\beta$ (see Appendix), for any bounded $\Omega\subset {\bf E}$
\begin{equation}
\begin{aligned}\label{rnon3}
\beta(S_{\bf A}(t)\Omega) \le \beta(S_{\bf A}(t){\bf P}_-\Omega) + \beta(S_{\bf A}(t){\bf Q}_+\Omega) = \beta(S_{\bf A}(t){\bf Q}_+\Omega),
\end{aligned}
\end{equation}
where the last inequality follows from the fact that the set $S_{\bf A}(t){\bf P}_-\Omega$ is relatively compact as a bounded subset of finite dimensional space ${\bf E}_0\oplus{\bf E}_-$. By \eqref{rnon1ab} and Lemma \ref{lem-noncom}, we deduce  that for any bounded $\Omega\subset {\bf E}_+$
\begin{equation*}
\beta_{{\bf E}_+} ({\bf Q}_+\Omega) = \beta ({\bf Q}_+\Omega) \text{ \ \ and \ \ } \beta_{{\bf E}_+}(S_{\bf A}(t){\bf Q}_+\Omega) = \beta(S_{\bf A}(t){\bf Q}_+\Omega) \ \ \text{ for } \ \ t\ge 0.
\end{equation*}
Therefore, combining this with \eqref{rnon1}, \eqref{rnon3} and \eqref{rnon1ab} yields
\begin{align*}
\beta(S_{\bf A}(t)\Omega) & \le \beta(S_{\bf A}(t){\bf Q}_+\Omega) = \beta_{{\bf E}_+}(S_{\bf A}(t){\bf Q}_+\Omega) \\
& \le e^{-\delta t} \beta_{{\bf E}_+} ({\bf Q}_+\Omega) = e^{-\delta t} \beta ({\bf Q}_+\Omega) \le e^{-\delta t} \beta (\Omega),
\end{align*}
which completes the proof of point $(b)$. \hfill $\square$ \\

In the following proposition we collect important facts concerning existence, continuity and compactness for the mild solutions of the equation \eqref{row-abstr}.

\begin{proposition}\label{th-exist1}
Under the above assumptions the following assertions hold.
\begin{itemize}
\item[(a)] For every $s\in[0,1]$ and $(x,y)\in {\bf E}$, the equation \eqref{row-abstr} admits a unique mild solution $w(\,\cdot\,;s, (x,y)):[0,+\infty)\to {\bf E}$ starting at $(x,y)$. Therefore, for any $t\ge 0$, one can define \emph{the Poincar\'e operator} associated with this equation as the map ${\bf \Phi}_t:[0,1]\times {\bf E} \to {\bf E}$ given by $${\bf \Phi}_t(s, x) := w(t; s, (x,y)) \text{ \ for \ } s\in [0,1] \ \text{ and } \ (x,y)\in {\bf E}.$$
\item[(b)] If sequences $(x_n,y_n)$ in ${\bf E}$ and $(s_n)$ in $[0,1]$ are such that $(x_n,y_n)\to (x_0,y_0)$ in ${\bf E}$ and $s_n\to s_0$ as $n\to +\infty$, then
\begin{equation*}
w(t;s_n, (x_n,y_n)) \to w(t;s_0, (x_0,y_0)) \qquad \mathrm{as} \quad n\to +\infty,
\end{equation*}
for $t\ge 0$, and this convergence in uniform on bounded subsets of $[0,+\infty)$. \\[-8pt]
\item[(c)] If $|\cdot|$ is the norm on ${\bf E}$ obtained in Proposition \ref{lem-con-sem} and $\beta$ is the Hausdorff measure of noncompactness associated with this norm, then
    $$\beta({\bf \Phi}_t([0,1]\times \Omega)) \le e^{-\delta t}\beta(\Omega)$$
    for $t\ge 0$ and any bounded $\Omega\subset{\bf E}$, where $\delta > 0$ is the constant from \eqref{nieee1}.
\end{itemize}
\end{proposition}
\noindent\textbf{Proof.} The proof of point $(a)$ is a consequence of \cite[Theorem 3.3.3]{MR610244} and
\cite[Corollary 3.3.5]{MR610244}. Points $(b)$ and $(c)$ are consequences of \cite[Proposition 4.1]{cw-kok1} and Proposition \ref{lem-con-sem}. \hfill $\square$

\section{Degree formula for periodic solutions}

In this section we intend to prove {\em the degree formula for periodic solutions}, which is the main result of this paper. We will consider the following differential equation
\begin{equation}\label{eq-hyp-tt}
    \ddot u(t) = -A u(t) - c A \dot u(t) + \lambda u(t) + F(t,u(t)), \qquad t\ge 0,
\end{equation}
where $c>0$ is a damping constant, $\lambda$ is a real number, the operator $A:X\supset D(A) \to X$ is defined on a separable Banach space $X$ with the norm $\|\cdot\|$ and $F\colon [0,+\infty)\times X^\alpha \to X$ is a continuous map. We are interested in the situation when the equation \eqref{eq-hyp-tt} is at {\em resonance at infinity}, which means that $$\mathrm{Ker}\,(\lambda I - A)\neq \{0\} \text{ \ and \ } F \text{ \ is a bounded map.}$$
Throughout this section we assume that conditions $(A1)-(A3)$, $(F1)$, $(F3)$ are satisfied and furthermore:\\[2pt]
\noindent\makebox[9mm][l]{$(F4)$}\parbox[t][][t]{118mm}{there is $m > 0$ such that $\|F(t,x)\|\le m$ for $t\in [0,+\infty)$ and $x\in X^\alpha$,} \\[5pt]
\noindent\makebox[9mm][l]{$(F5)$}\parbox[t][][t]{118mm}{there is $T > 0$ such that $F(t + T,x) = F(t,x)$ for $t\in [0,+\infty)$ and $x\in X^\alpha$.} \\[5pt]
Let us write the equation \eqref{eq-hyp-tt} in the following form
\begin{equation}\label{diff:2kk22}
    \dot w(t) = -{\bf A} w(t) + {\bf F}(t, w(t)), \qquad t > 0
\end{equation}
where ${\bf A}:{\bf E}\supset D({\bf A})\to{\bf E}$ is a linear operator given by the formula \eqref{op-hyp} and ${\bf F}:[0,+\infty)\times {\bf E}\to{\bf E}$ is a map defined by
\begin{equation*}
    {\bf F}(t, (x,y)):=(0,F(t, x)) \ \ \text{for} \ t\in[0,+\infty), \ (x,y)\in{\bf E}.
\end{equation*}
Throughout this section the space ${\bf E} = X^\alpha\times X$ is equipped with the equivalent norm $|\cdot|$ obtained in Proposition \ref{lem-con-sem}. By Proposition \ref{th-exist1} $(a)$, for any $(x,y)\in{\bf E}$, there is a unique mild solution $w(\,\cdot\,;(x,y)):[0,+\infty)\to {\bf E}$ of \eqref{diff:2kk22} starting at $(x,y)$.
Given $T \ge 0$, let us define the Poincar\'e operator ${\bf \Phi}_T:{\bf E}\to{\bf E}$ given by
\begin{equation*}
    {\bf \Phi}_T(x,y):= w(T;(x,y)) \qquad \mathrm{for} \quad (x,y)\in{\bf E}.
\end{equation*}
Proposition \ref{th-exist1} $(b)$ and $(c)$ say that ${\bf \Phi}_T$ is continuous and
\begin{align*}
\beta({\bf \Phi}_T(\Omega)) \le e^{-\delta T}\beta(\Omega),
\end{align*}
for any bounded $\Omega\subset{\bf E}$, where $\beta$ is a Hausdorff measure of noncompactness associated with the norm $|\cdot|$ and the constant $\delta > 0$ is from \eqref{nieee1}.
\begin{remark}\label{rem-non-ex4}
If the equation \eqref{eq-hyp-tt} is at resonance at infinity, then there is a nonlinearity $F$ such that the equation
does not admit $T$-periodic mild solutions. \\[5pt]
To see this, let us take $F(t,x) = y_0$ for $t\in[0,+\infty)$ and $x\in {\bf E}$, where $y_0\in\mathrm{Ker}\,(\lambda I - A)\setminus\{0\}$. If $w:[0,+\infty)\to {\bf E}$ is a $T$-periodic solution for \eqref{diff:2kk22}, then it satisfies the integral formula
\begin{equation}\label{eq111}
w(t) = S_{\bf A}(t)w(0) + \int_0^t S_{\bf A}(t - \tau)(0,y_0)\,d\tau\qquad \mathrm{for} \quad t\ge 0.
\end{equation}
Consider the direct sum decomposition ${\bf E}:={\bf E}_-\oplus{\bf E}_0\oplus{\bf E}_+$ obtained in Theorem \ref{th-spec-dec} together with the corresponding continuous projections ${\bf P}$, ${\bf Q}_+$, ${\bf Q}_-$. Let us equip ${\bf E}_0=\mathrm{Ker}\,(\lambda I - A)\times \mathrm{Ker}\,(\lambda I - A)$ with the following norm
\begin{equation}\label{wz-norm}
    \|(x,y)\|_{{\bf E}_0} := \|x\|_H + \|y\|_H \qquad \mathrm{for} \quad (x,y)\in{\bf E}_0.
\end{equation}
Acting on the equation \eqref{eq111} by ${\bf P}$ and using \eqref{ds1-hyp}, \eqref{dss1}, \eqref{proj} one has $${\bf P} w(t) = S_{{\bf A}_0}(t){\bf P} w(0) + \int_0^t S_{{\bf A}_0}(t - \tau)(0, y_0)\,d\tau\qquad \mathrm{for} \quad t\ge 0.$$
Let $(u_0, v_0):\mathbb{R}\to {\bf E}_0$, be a map given by $(u_0(t), v_0(t)) := {\bf P} w(t)$ for $t\ge 0$. Since ${\bf A}_0$ is a bounded operator on a finite dimensional space ${\bf E}_0$, it follows that $(u_0, v_0)$ is of class $C^1$ and
\begin{equation*}
\left\{ \begin{aligned}
\dot u_0(t) & = v_0(t), & \qquad & t \ge 0, \\
\dot v_0(t) & = - c \lambda v_0(t) + y_0,  &\qquad& t \ge 0.
\end{aligned}\right.
\end{equation*}
As a consequence, for any $t\ge 0$, one has
\begin{align*}
\frac{d}{dt} (\<u_0(t), c \lambda y_0\>_H + \<v_0(t), y_0\>_H) = \<v_0(t), c \lambda y_0\>_H + \<- c \lambda v_0(t) + y_0, y_0\>_H = \|y_0\|_H^2.
\end{align*}
Therefore, if $w:[0,+\infty)\to {\bf E}$ is a $T$-periodic solution of \eqref{diff:2kk22}, then
\begin{equation*}
   0 = \<u_0(T), c \lambda y_0\>_H + \<v_0(T), y_0\>_H - \<u_0(0), c \lambda y_0\>_H - \<v_0(0), y_0\>_H = T\|y_0\|_H^2,
\end{equation*}
which is impossible, because $y_0 \neq 0$. \hfill $\square$
\end{remark}

To overcome difficulties presented in Remark \ref{rem-non-ex4} we impose additional assumptions on nonlinearity $F$. Recalling that $X^\alpha_+$, $X^\alpha_-$ and $X_0$ are subspaces from Remark \ref{rem-proj}, we equip the spaces $X^\alpha_+\oplus X^\alpha_-$ and $X_0$ with the norms $\|\cdot\|_\alpha$ and $\|\cdot\|_H$, respectively. Then
we introduce the following geometric conditions:
\begin{equation*}\leqno{(G1)}
\ \left\{\begin{aligned}
& \text{for any balls } B_1\subset X^\alpha_+\oplus X^\alpha_- \text{ and } B_2\subset X_0 \text{ there is } R > 0 \text{ such that } \\
& \<F(t, x + y), x\>_H > - \<F(t, x + y), z\>_H \\
& \text{for } (t,y,z)\in [0,T]\times B_1 \times B_2, \ \text{and} \ x\in X_0 \text{ with } \|x\|_H\ge R,
\end{aligned}\right.
\end{equation*}
\begin{equation*}\leqno{(G2)}
\ \left\{\begin{aligned}
& \text{for any balls } B_1\subset X^\alpha_+\oplus X^\alpha_- \text{ and } B_2\subset X_0 \text{ there is } R > 0 \text{ such that } \\
& \<F(t, x + y), x\>_H < - \<F(t,x + y), z\>_H \\
& \text{for } (t,y,z)\in [0,T]\times B_1 \times B_2, \ \text{and} \ x\in  X_0 \text{ with } \|x\|_H\ge R.
\end{aligned}\right.
\end{equation*}
Now we are ready to prove the following \emph{index formula for periodic solutions}, which uses the conditions $(G1)$ and $(G2)$ to determine the topological degree for the condensing vector field $I - {\bf \Phi}_T$ (see Appendix). We recall that by Remark \ref{rem-pom} the spectrum $\sigma(A)$ consists of sequence of positive eigenvalues $(\lambda_i)_{i\ge 1}$ and $\dim\mathrm{Ker}\,(\lambda_i I - A) < +\infty$ for $i\ge 1$.
\begin{theorem}\label{th-reso-m-h}
Let $\lambda = \lambda_k$ be the $k$-th eigenvalue of the operator $A$.
\begin{enumerate}
\item[(i)] If condition $(G1)$ is satisfied, then there is an open bounded set $W \subset {\bf E}$ such that ${\bf \Phi}_T(x,y)\neq (x,y)$ for $(x,y)\in\partial W$ and
\begin{equation*}
    \mathrm{deg_{C}}(I - {\bf \Phi}_T, W) = (-1)^{d_k}.
\end{equation*}
\item[(ii)] If condition $(G2)$ is satisfied, then there is an open bounded set $W \subset {\bf E}$ such that ${\bf \Phi}_T(x,y)\neq (x,y)$ for $(x,y)\in\partial W$ and
\begin{equation*}
    \mathrm{deg_{C}}(I - {\bf \Phi}_T, W) = (-1)^{d_{k-1}}.
\end{equation*}
\end{enumerate}
Here $d_l:=\sum_{i=1}^l \dim\mathrm{Ker}\,(\lambda_i I - A)$ for $l\ge 1$ with the exceptional case $d_0:= 0$.
\end{theorem}

\subsection{Preparation to the proof of Theorem \ref{th-reso-m-h}}

In the proof we will use the following family of differential equations
\begin{equation}\label{diff:2kkaabb}
\dot w(t) = - {\bf A} w(t) + {\bf G}(s,t,w(t)), \qquad t > 0
\end{equation}
where ${\bf G}:[0,1]\times [0,+\infty)\times {\bf E} \to {\bf E}$ is a map given by
$${\bf G}(s,t,(x,y)):= (0,G(s,t,x)) := (0, PF(t,s Q x + P x) + s QF(t,s Q x + P x))$$
for $s\in[0,1]$, $t\in [0,+\infty)$ and $(x,y)\in{\bf E}$.
\begin{remark}\label{rem-com11}
It is not difficult to see that $G$ satisfies assumption $(F1)-(F3)$. Furthermore there is a constant $m_0 > 0$ such that
\begin{equation}\label{roow}
\|{\bf G}(s,t, z)\|_{\bf E} \le m_0 \qquad \mathrm{for} \quad s\in[0,1], \ t\in[0,+\infty), \ z\in{\bf E},
\end{equation}
which is a consequence of assumption $(F4)$. \hfill $\square$
\end{remark}
In the view of the above remark, Proposition \ref{th-exist1} $(a)$ asserts that for aby $s\in [0,1]$ and $(x,y)\in{\bf E}$ there is a map $w(\,\cdot\,;s, (x,y)):[0,+\infty)\to {\bf E}$ which ia a mild solution of \eqref{diff:2kkaabb} starting at $(x,y)$. Let ${\bf \Psi}_T:[0,1]\times{\bf E}\to{\bf E}$ be the associated Poincar\'e operator given by
\begin{equation*}
    {\bf \Psi}_T(s,(x,y)):= w(T;s,(x,y)) \qquad \mathrm{for} \quad s\in [0,1], \ (x,y)\in{\bf E}.
\end{equation*}
Remark \ref{rem-com11} and Proposition \ref{th-exist1} $(b)$ and $(c)$ imply that ${\bf \Psi}_T$ is continuous and
\begin{align*}
\beta({\bf \Psi}_T([0,1]\times \Omega)) \le e^{-\delta T}\beta(\Omega)
\end{align*}
for any bounded set $\Omega\subset{\bf E}$ where the constant $\delta$ is from \eqref{nieee1}. The following lemma provides some {\em a priori} bounds on periodic solutions for the equation \eqref{diff:2kkaabb}.

\begin{lemma}\label{lem-bounkk}
There is a constant $R > 0$ such that, if $w:[0,+\infty)\to {\bf E}$, is a $T$-periodic mild solution of \eqref{diff:2kkaabb}, then
\begin{equation*}
    \|{\bf Q} w(t)\|_{\bf E} \le R \qquad \mathrm{for} \quad t\in [0,T],
\end{equation*}
where ${\bf Q} = {\bf Q}_+ + {\bf Q}_-$ and ${\bf Q}_+$, ${\bf Q}_-$ are projections obtained in Theorem \ref{th-spec-dec}.
\end{lemma}
\noindent\textbf{Proof.} Let us first observe that $w(t) = w(t + kT)$ for any integer $k \ge 0$ and $t\in [0,+\infty)$, which follows from the fact that $w$ is $T$-periodic and $G$ satisfies $(F5)$. This in turn implies that
\begin{equation}\label{rere}
w(t) = S_{\bf A}(kT)w(t) + \int_t^{t + kT} S_{\bf A}(t + kT - \tau){\bf G}(s,\tau, w(\tau)) \, d\tau
\end{equation}
for $t\ge 0$ and $n\ge 1$. Acting on this equation by ${\bf Q}_+$ and using \eqref{ds1-hyp}, we infer that
\begin{equation*}
{\bf Q}_+w(t) = S_{\bf A}(kT){\bf Q}_+ w(t) + \int_t^{t + kT} S_{\bf A}(t + kT - \tau) {\bf Q}_+ {\bf G}(s,\tau, w(\tau)) \, d \tau
\end{equation*}
for $t\ge 0$ and $n\ge 1$. Therefore, by \eqref{nieee1}, one has
\begin{equation*}
\begin{aligned}
\|{\bf Q}_+w(t)\|_{\bf E} & \le \|S_{\bf A}(kT){\bf Q}_+w(t)\|_{\bf E} \\
& \hspace{5mm} + \int_t^{t + k T} \|S_{\bf A}(t + k T - \tau){\bf Q}_+{\bf G}(s,\tau,w(\tau))\|_{\bf E} \, d \tau \\
& \le \|S_{\bf A}(k T){\bf Q}_+w(t)\|_{\bf E} + \int_t^{t + k T} M e^{- \delta (t + k T - \tau)} \|{\bf Q}_+{\bf G}(s,\tau,w(\tau))\|_{\bf E} \, d \tau \\
& \le \|S_{\bf A}(k T){\bf Q}_+w(t)\|_{\bf E} + \int_t^{t + k T} m_0 M  \|{\bf Q}_+\|_{L({\bf E})} e^{- \delta (t + k T - \tau)} \, d \tau \\
& \le M e^{- \delta kT } \, \|{\bf Q}_+w(t)\|_{\bf E} + \int_t^{t + k T} m_0M  \|{\bf Q}_+\|_{L({\bf E})}e^{- \delta (t + k T - \tau)}\, d \tau \\
& \le M e^{- \delta kT } \, \|{\bf Q}_+w(t)\|_{\bf E} + m_0 M \|{\bf Q}_+\|_{L({\bf E})}\left(1 - e^{- \delta k T}\right)/\delta,
\end{aligned}
\end{equation*}
where $m_0 > 0$ is a constant from \eqref{roow}. In a consequence, for any $t\in [0,T]$ and integer $k > 0$, we obtain
\begin{equation*}
\|{\bf Q}_+w(t)\|_{\bf E} \le M e^{- \delta kT } \, \|{\bf Q}_+w(t)\|_{\bf E} + m_0 M \|{\bf Q}_+\|_{L({\bf E})}\left(1 - e^{- \delta k T}\right)/\delta
\end{equation*}
and hence, letting $k\to +\infty$ yields
\begin{equation}\label{eq:31kk}
\|{\bf Q}_+w(t)\|_{\bf E} \le m_0 M \|{\bf Q}_+\|_{L({\bf E})}/\delta :=R_1\qquad \mathrm{for} \quad t\in [0,T].
\end{equation}
On the other hand, acting on the equation \eqref{rere} by ${\bf Q}_-$ and applying \eqref{ds1-hyp}, gives
\begin{equation*}
{\bf Q}_-w(t) = S_{\bf A}(k T){\bf Q}_-w(t) + \int_t^{t + kT} S_{\bf A}(t + kT - \tau){\bf Q}_-{\bf G}(s,\tau, w(\tau)) \, d \tau.
\end{equation*}
By Corollary \ref{cor-dec-2aba} $(ii)$, the semigroup $\{S_{\bf A}(t)\}_{t\ge 0}$ extends on ${\bf E}_-$ to a $C_0$ group of bounded operators. Hence, for any $t\in [0,T]$ and integer $k\ge 1$, one has
\begin{equation*}
S_{\bf A}(-k T){\bf Q}_-w(t) = {\bf Q}_-w(t) + \int_t^{t + kT} S_{\bf A}(t - \tau){\bf Q}_-{\bf G}(s,\tau, w(\tau)) \, d \tau,
\end{equation*}
which together with \eqref{nieee2} and \eqref{roow} gives
\begin{align*}
\|{\bf Q}_-w(t)\|_{\bf E} & \le \|S_{\bf A}(- k T) {\bf Q}_-w(t)\|_{\bf E} \\
& \qquad + \int_t^{t + k T} \|S_{\bf A}(t - \tau){\bf Q}_-{\bf G}(s,\tau,w(\tau))\| \, d \tau \\
& \le M  \, e^{-\delta k T} \|{\bf Q}_-w(t)\|_{\bf E} + M\int_t^{t + kT} e^{\delta (t - \tau)}\|{\bf Q}_-{\bf G}(s,\tau,w(\tau))\| \, d \tau \\
& \le M  \, e^{-\delta k T} \|{\bf Q}_-w(t)\|_{\bf E} + m_0 M \int_t^{t + kT} \|{\bf Q}_-\|_{L({\bf E})} e^{\delta (t - \tau)} \, d \tau \\
& = M  \, e^{-\delta k T} \|{\bf Q}_-w(t)\|_{\bf E} + m_0 M \|{\bf Q}_-\|_{L({\bf E})}\left(1 - e^{ - \delta k T}\right)/\delta.
\end{align*}
Hence, passing to the limit with $k\to +\infty$, we obtain
\begin{equation}\label{eq:32kk}
\|{\bf Q}_-w(t)\|_{\bf E} \le m_0 M \|{\bf Q}_-\|_{L({\bf E})}/\delta :=R_2 \qquad \mathrm{for} \quad t\in [0,T].
\end{equation}
Finally, by \eqref{eq:31kk} and \eqref{eq:32kk}, we find that
$$\|{\bf Q} w(t)\|_{\bf E} \le \|{\bf Q}_+ w(t)\|_{\bf E} + \|{\bf Q}_-w(t)\|_{\bf E} \le R_1 + R_2:=R,$$ for $t\in[0,T]$, and the proof is completed. \hfill $\square$ \\

We proceed to consider the following version  of the equation \eqref{diff:2kkaabb}:
\begin{equation}\label{diff:2kkaabbcc}
\dot w(t) = - \mu{\bf A} w(t) + \mu{\bf G}(s,t,w(t)), \qquad t > 0
\end{equation}
where $\mu\in(0,1]$ is a parameter. We will need the following lemma.

\begin{lemma}\label{lem-per-bound11}
There is a constant $R > 0$ with the property that, if $w=(u,v):[0,+\infty)\to {\bf E}$ is a $T$-periodic mild solution of \eqref{diff:2kkaabbcc} for some $\mu\in(0,1]$, then $\|P v(t)\|_H\le R$ for $t\in[0,T]$.
\end{lemma}
\noindent\textbf{Proof.} Acting by the operator ${\bf P}$ on the equation
\begin{equation*}
w(t) = S_{\bf A}(\mu t)w(0) + \int_0^t S_{\bf A}(\mu(t - \tau))(0, \mu G(s, t,u(\tau)))\,d\tau\qquad \mathrm{for} \quad t\ge 0.
\end{equation*}
and using \eqref{ds1-hyp}, \eqref{dss1}, \eqref{proj} one has
$${\bf P} w(t) = S_{{\bf A}_0}(\mu t){\bf P} w(0) + \int_0^t S_{{\bf A}_0}(\mu(t - \tau))(0, \mu PF(t,s Q u(\tau) + P u(\tau)))\,d\tau \text{ \ for \ } t\ge 0.$$
Since ${\bf A}_0$ is a bounded operator on a finite dimensional space ${\bf E}_0$, it follows that the map $(u_0, v_0) := {\bf P} w = (Pu, Pv)$ continuously differentiable on $\mathbb{R}$ and
\begin{equation}
\left\{\begin{aligned}\label{ror1}
\dot u_0(t) & = \mu v_0(t) \qquad t \ge 0, \\
\dot v_0(t) & =  - c\mu\lambda v_0(t) + \mu PF(t,s Q u(t) + P u(t)) \qquad t \ge 0.
\end{aligned} \right.
\end{equation}
Let $k:[0,+\infty)\to X_0$, where we recall that $X_0 = \mathrm{Ker}\,(\lambda I - A)$, be a map given by $$k(t):= PF(t,s Q u(t) + P u(t)) \qquad \mathrm{for} \quad t\ge 0.$$ Since $X_0$ is a finite dimensional space, we have the norm equivalence $$c \|x\|_H \le \|x\| \le C\|x\|_H \qquad \mathrm{for} \quad x\in X_0,$$ where $c,C > 0$ are constants. This leads us to
\begin{equation}
\begin{aligned}\label{eeqqq}
    \|k(t)\|_H & = \|PF(t,s Q u(t) + P u(t))\|_H \\
    & \le C\|PF(t,s Q u(t) + P u(t))\| \le C\|P\| m := m_1 \quad \text{ for \ } t\ge 0
\end{aligned}
\end{equation}
where the last inequality is a consequence of $(F4)$. Chose $R_0 > 0$ such that
\begin{equation}\label{njnj2}
-c\lambda R^2 + m_1 R < 0 \qquad \mathrm{for} \quad R\ge R_0
\end{equation}
and observe that, by \eqref{ror1} and \eqref{eeqqq}, one has
\begin{equation}
\begin{aligned}\label{roro2}
\frac{d}{dt}\frac{1}{2}\|v_0(t)\|_H^2 & = \<\dot v_0(t), v_0(t)\>_H=\<- c\mu\lambda v_0(t)+\mu k(t),v_0(t)\>_H \\
& = - c\mu\lambda \|v_0(t)\|_H^2 + \mu \<k(t), v_0(t)\>_H \\
& \le - c\mu\lambda \|v_0(t)\|_H^2 + \mu m_1\|v_0(t)\|_H  \text{ \ \ for \ } t\ge 0.
\end{aligned}
\end{equation}
Let us observe that there is $t_0 \ge 0$ such that $\|v_0(t_0)\|_H\le R_0+1$, because otherwise, we would have $\|v_0(t)\|_H \ge R_0 + 1$ for $t\ge 0$, which together with \eqref{roro2} and \eqref{njnj2} would give
$$\frac{d}{dt}\frac{1}{2}\|v_0(t)\|_H^2  < 0 \qquad \mathrm{for} \quad t\ge 0.$$ This is impossible because $v_0$ is $T$-periodic.
Let us now define
\begin{align*}
& D(0,R_0+1) := \{x\in \mathrm{Ker}\, (\lambda I - A) \ | \ \|x\|_H \le R_0+1\}, \\
& A :=\{\delta \ge 0 \ | \ v_0([t_0,t_0+\delta]) \subset D(0,R_0+1)\}.
\end{align*}
Then $A$ is a nonempty set because $0\in A$ and hence one can consider $s:=\sup A$. If $s < +\infty$ then $s\in A$ because $D(0,R_0+1)$ is a closed set. Furthermore, by continuity of $v_0$, one has $\|v_0(t_0+s)\|_H = R_0+1$. Let us note that the equation \eqref{roro2} implies
\begin{equation*}
\begin{aligned}
\frac{d}{dt}\frac{1}{2}\|v_0(t)\|_H^2|_{t = t_0 + s}  \le - c\mu\lambda \|v_0(t_0 + s)\|_H^2 + \mu m_1 \|v_0(t_0 + s)\|_H < 0
\end{aligned}
\end{equation*}
and consequently there is $\delta_0 > 0$ such that $v_0([t_0,t_0+ s+ \delta_0])\subset D(0,R_0+1)$, which is impossible in view of the definition of $s$. Hence $s=\sup A = +\infty$ which together with $T$-periodicity of $v_0$ gives
$$\|Pv(t)\|_H = \|v_0(t)\|_H \le R:= R_0 + 1 \qquad \mathrm{for} \quad t\in [0,+\infty)$$ and the proof of lemma is completed. \hfill $\square$ \\[5pt]

The following lemma provides another {\em a priori} bounds on periodic solutions for the equation \eqref{diff:2kkaabbcc}. In the proof of this lemma for the first time we will use the geometrical conditions $(G1)$ and $(G2)$.

\begin{lemma}\label{lem-per-bound22}
There is $R > 0$ such that for any $\mu\in(0,1]$, the equation \eqref{diff:2kkaabbcc} does not admit $T$-periodic mild solutions $w:[0,+\infty)\to {\bf E}$ with $\|{\bf P} w(0)\|_{{\bf E}_0} \ge R$.
\end{lemma}
\noindent\textbf{Proof.} Suppose contrary to our claim that there are sequences $(\mu_n)$ in $(0,1]$, $(s_n)$ in $[0,1]$ and a sequence of $T$-periodic mild solutions $w_n=(u_n,v_n):[0,+\infty)\to {\bf E}$ of \eqref{diff:2kkaabbcc} such that
$\|{\bf P} w_n(0)\|_{{\bf E}_0} \to +\infty$ as $n\to+\infty$. Acting by ${\bf P}$ on the equation
\begin{equation*}
w_n(t) = S_{\mu_n{\bf A}}(t)w_n(0) + \mu_n\int_0^t S_{\mu_n{\bf A}}(\mu_n(t - \tau))(0, G(s_n, \tau ,u_n(\tau)))\,d\tau\qquad \mathrm{for} \quad t\ge 0.
\end{equation*}
and using \eqref{ds1-hyp}, \eqref{dss1}, \eqref{proj} one obtains
$${\bf P} w_n(t) = S_{\mu_n{\bf A}_0}(t){\bf P} w_n(0) + \mu_n\int_0^t S_{\mu_n{\bf A}_0}(t - \tau)(0, PF(\tau,s_n Q u_n(\tau) + P u_n(\tau)))\,d\tau.$$
Without loss of generality we can assume also that $\mu_n\to \mu_0$ and $s_n\to s_0$ as $n\to +\infty$, where $\mu_0,s_0\in [0,1]$. Let us define $z_n := {\bf P} w_n(0)/\|{\bf P} w_n(0)\|_{{\bf E}_0}$ for $n\ge 1$. Since $(z_n)_{n\ge 1}$ is a bounded sequence contained in the finite dimensional space ${\bf E}_0$, we can also assume that there is $z_0\in{\bf E}_0$ such that
$$z_n = {\bf P} w_n(0)/\|{\bf P} w_n(0)\|_{{\bf E}_0} \to z_0 \ \ \text{ as }\ n\to +\infty.$$
Since the maps $w_n$ are $T$-periodic, assumption $(F5)$ implies that $w_n(t) = w_n(t + T)$ for $t\ge 0$, which in consequence gives $w_n(0) = w_n(kT)$ for any integer $k\ge 1$.
Let us note that we can choose a sequence of integers $(k_n)_{n\ge 1}$ such that $k_n\mu_n T\to t_0 > 0$ as $n\to +\infty$. Indeed, if $\mu_0\neq 0$ then it is enough to take $k_n = 1$. If $\mu_0 = 0$ then we define $k_n:= \lfloor t_0/(\mu_n T)\rfloor$, where $\lfloor x\rfloor$ is the integer part of $x$. Therefore
\begin{equation}\label{row-gr}
z_n = S_{{\bf A}_0}(\mu_n k_nT) z_n + y_n(k_nT) \qquad \mathrm{for} \quad n\ge 1.
\end{equation}
where, for $n\ge 1$, we denote
$$y_n(t):=\mu_n \int_0^tS_{\mu_n{\bf A}_0}(t-\tau)(0, PF(\tau,s_n Q u_n(\tau) + P u_n(\tau)))/\|w_n(0)\|_{{\bf E}_0}\,d\tau.$$
Since $\{S_{{\bf A}_0}(t)\}_{t\ge 1}$ is a $C_0$ semigroup, there are $\omega\in\mathbb{R}$ and $M > 0$ such that
$$\|S_{{\bf A}_0}(t)\|_{{\bf E}_0} \le M e^{\omega t} \qquad \mathrm{for} \quad t\ge 0.$$
Using this together with \eqref{eeqqq}, we infer that
\begin{equation*}
\hspace{-2mm}\begin{aligned}
\|y_n(t)\|_{{\bf E}_0} & \le \mu_n\int_0^{t} \|S_{\mu_n{\bf A}_0}(t - \tau)(0, PF(\tau,s_n Q u_n(\tau) + P u_n(\tau)))\|_{{\bf E}_0}/\|w_n(0)\|_{{\bf E}_0} \, d \tau \\ & \le \mu_n\int_0^{t} M e^{\omega\mu_n(t - \tau)} \|(0, PF(\tau,s_n Q u_n(\tau) + P u_n(\tau)))\|_{{\bf E}_0}/\|w_n(0)\|_{{\bf E}_0} \, d \tau \\
& \le \mu_n\int_0^{t} M e^{\omega\mu_n(t - \tau)} m_1 /\|w_n(0)\|_{{\bf E}_0} \, d \tau \qquad \mathrm{for} \quad n\ge 1, \ t\ge 0,
\end{aligned}
\end{equation*}
which after integration of the last term gives
\begin{equation}\label{zbie-jedn}
\|y_n(t)\|_{{\bf E}_0} \le \frac{M m_1}{\omega \|w_n(0)\|_{{\bf E}_0}}(e^{\omega\mu_n t} - 1) \qquad \mathrm{for} \quad n\ge 1, \ t\ge 0.
\end{equation}
In particular, this implies that $y_n(k_n T)\to 0$ as $n\to +\infty$. Letting $n\to +\infty$ in the equation \eqref{row-gr}, yields
$$z_0 = S_{{\bf A}_0}(t_0) z_0 = S_{\bf A}(t_0) z_0 \qquad\text{for some}\quad t_0 > 0,$$
where the last equality follows from \eqref{dss1}. Hence, by Theorem \ref{th-spec} $(iii)$, one has $$z_0\in\mathrm{Ker}\,{\bf A} = \{(x,0) \ | \ x\in \mathrm{Ker}\,(\lambda I -A)\}$$ and consequently $z_0 = (z^1_0, 0)$, where $z^1_0\in \mathrm{Ker}\,(\lambda I -A)$. Furthermore $z_0\in{\bf E}_0$ and $$z_0 = S_{\bf A} (t) z_0 = S_{{\bf A}_0}(t) z_0 \qquad \mathrm{for} \quad t > 0$$ because $z_0\in\mathrm{Ker}\,{\bf A}$. Combining this with \eqref{zbie-jedn} yields $$(P u_n(t), P v_n(t))/\|{\bf P} w_n(0)\|_{{\bf E}_0} = {\bf P} w_n(t)/\|{\bf P} w_n(0)\|_{{\bf E}_0} \to z_0 = (z^1_0, 0) \text{ \ as \ } n\to +\infty,$$ uniformly for $t\in [0,T]$, which by \eqref{wz-norm} implies that
\begin{equation}\label{rowro1}
    P u_n(t)/\|{\bf P} w_n(0)\|_{{\bf E}_0} \to z^1_0 \ \ \text{ uniformly for } \ t\in [0,T] \ \text{ \ in the norm  \ } \|\cdot\|_H.
\end{equation}
Since $\|z_n\|_{{\bf E}_0} = 1$ for $n\ge 1$, it follows that $\|z^1_0\|_H = \|z_0\|_{{\bf E}_0} = 1$. As the operator ${\bf A}_0$ is bounded, the maps ${\bf P} w_n=(P u_n, P v_n):[0,+\infty)\to {\bf E}_0$ are of class $C^1$ and
\begin{equation*}
\left\{\begin{aligned}
\frac{d}{dt} P u_n(t) & = \mu_n P v_n(t), && t \ge 0 \\
\frac{d}{dt} P v_n(t) & =  - c\mu_n \lambda P v_n(t) + \mu_n PF(t,s_n Q u_n(t) + P u_n(t)), && t \ge 0.
\end{aligned} \right.
\end{equation*}
This implies that, for any $n\ge 1$ and $t\ge 0$
\begin{align*}
\frac{d}{dt}\frac{1}{2} \left\|c\lambda P u_n(t) + P v_n(t)\right\|_H^2
& = \<c\lambda \frac{d}{dt} P u_n(t) + \frac{d}{dt} P v_n(t), c\lambda P u_n(t) + P v_n(t)\>_H \\
& = \mu_n \<PF(t, s_n Q u_n(t) + P u_n(t)), c\lambda P u_n(t) + P v_n(t)\>_H,
\end{align*}
which after integration gives
\begin{align*}
& \frac{1}{2} \left\|c\lambda P u_n(T) + P v_n(T)\right\|_H^2 - \left\|c\lambda P u_n(0) + P v_n(0)\right\|_H^2 \\
& \hspace{30mm} = \mu_n \int_0^T\<PF(\tau, s_n Q u_n(\tau) + P u_n(\tau)), c\lambda P u_n(\tau) + P v_n(\tau)\>_H \, d \tau \end{align*}
Since the maps $(u_n,v_n)$ are $T$-periodic, we infer that
\begin{align}\label{row-cal}
\int_0^T\<PF(\tau, s_n Q u_n(\tau) + P u_n(\tau)), c\lambda P u_n(\tau) + P v_n(\tau)\>_H \, d \tau = 0.
\end{align}
By Lemma \ref{lem-bounkk} there is $R_1 > 0$ such that
$$\|s_n Q u_n(t)\|_\alpha \le \|{\bf Q} w_n(t)\|_{\bf E} \le R_1 \qquad \mathrm{for} \quad n\ge 1, \ t\in [0,+\infty).$$ On the other hand, by Lemma \ref{lem-per-bound11}, one can choose $R_2 > 0$ such that $\|P v_n(t)\|_H\le R_2$ for $t\in[0,T]$.
Let us define
\begin{align*}
B_1:=\{y\in X^\alpha_+\oplus X^\alpha_- \ | \ \|y\|_\alpha < R_1\} \text{ \ and \ }
B_2:=\{y\in X_0\ | \ \|y\|_H < R_2 /(c\lambda) \}.
\end{align*}
Using geometrical conditions $(G1)$, $(G2)$ and Proposition \ref{th:10} $(iii)$, one can take a constant $R > 0$ such that
$$\<PF(t, x + y), x\>_H > - \<PF(t,x + y), z\>_H$$ for $(t,y, z)\in [0,T]\times B_1\times B_2$ and $x\in X_0$ with $\|x\|_H\ge R$ if $(G1)$ is satisfied and
$$\<PF(t, x + y), x\>_H < - \<PF(t, x + y), z\>_H$$ for $(t,y, z)\in [0,T]\times B_1\times B_2$ and $x\in X_0$ with $\|x\|_H\ge R$ if $(G2)$ is satisfied. Since $\|z^1_0\|_H = 1$, by the convergence \eqref{rowro1}, we can choose $n_0 \ge 1$ such that
$$\left\|\frac{P u_n(t)}{\|{\bf P} w_n(0)\|_{{\bf E}_0}}\right\|_H \ge \|z^1_0\|_H - 1/2 = 1/2 \qquad \mathrm{for} \quad n\ge n_0 \text{ \ and \ }t\in[0,T].$$ Using this together with the fact that $\|{\bf P} w_n(0)\|_{{\bf E}_0} \to +\infty$ as $n\to +\infty$, we can  increase $n_0 \ge 1$ if necessary, and obtain
$$\|P u_n(t)\|_H = \|{\bf P} w_n(0)\|_{{\bf E}_0} \cdot \left\|\frac{P u_n(t)}{\|{\bf P} w_n(0)\|_{{\bf E}_0}}\right\|_H \ge \frac{1}{2} \|{\bf P} w_n(0)\|_{{\bf E}_0} \ge R$$ for $n\ge n_0$ and $t\in[0,T]$. This yields
\begin{align*}
\int_0^T \< PF(\tau, s_n Q u_n(\tau) + P u_n(\tau)),  c\lambda P u_n(\tau) + P v_n(\tau)\>_H \, d \tau  > 0 \text{ \ for \ } n\ge n_0,
\end{align*}
if condition $(G1)$ is satisfied, and
\begin{align*}
\int_0^T \< PF(\tau, s_n Q u_n(\tau) + P u_n(\tau)),  c\lambda P u_n(\tau) + P v_n(\tau)\>_H \, d \tau  < 0 \text{ \ for \ } n\ge n_0,
\end{align*}
if condition $(G2)$ is satisfied. Each of the above inequalities contradicts \eqref{row-cal} and thus the proof of lemma is completed. \hfill $\square$ \\

Let us now define {\em the averaging map} $\widehat{{\bf F}}:{\bf E}_0\to{\bf E}_0$ by the formula
\begin{equation*}
\widehat{{\bf F}}(x,y):= (0,\widehat{F}(x)) \qquad \mathrm{for} \quad (x,y)\in {\bf E}_0 = X_0 \times X_0
\end{equation*}
where we remind that in the above expression $X_0 = \mathrm{Ker}\,(\lambda I - A)$ and furthermore
\begin{equation*}
\widehat{F}(x) := \int_0^T PF(\tau,x)\,d \tau \qquad\mathrm{for}\quad x\in X_0,
\end{equation*}

In the following lemma we describe relation between the Brouwer degree of the averaging map $- {\bf A}_0 + \widehat{{\bf F}}$ and the geometrical conditions $(G1)$ and $(G2)$.
\begin{lemma}\label{th-guid-fun-kk22}
There is a constant $R_0 > 0$ such that the following assertions hold.
\begin{enumerate}
\item[(i)] If condition $(G1)$ is satisfied, then $-{\bf A}_0 z + \widehat{{\bf F}}(z) \neq 0$ for $\|z\|_{{\bf E}_0} \ge R_0$ and
\begin{equation*}
    \mathrm{deg_B}({\bf A}_0 - \widehat{{\bf F}}, U_R) = (-1)^{\dim X_0} \quad\mathrm{for}\quad R\ge R_0.
\end{equation*}
\item[(ii)] If condition $(G2)$ is satisfied, then $-{\bf A}_0 z + \widehat{{\bf F}}(z) \neq 0$ for $\|z\|_{{\bf E}_0} \ge R_0$ and
\begin{equation*}
    \mathrm{deg_B}({\bf A}_0 - \widehat{{\bf F}}, U_R) = 1 \quad\mathrm{for}\quad R\ge R_0.
\end{equation*}
\end{enumerate}
Here we denote $U_r := \{(x,y)\in {\bf E}_0 \ | \ \|x\|_H + \|y\|_H < r\}$ for $r>0$.
\end{lemma}
Before we start the proof we invoke the following lemma. For the proof we refer the reader to \cite[Theorem 5.2]{Kok2}.
\begin{lemma}\label{th-guid-fun-kk}
There is a constant $R_0 > 0$ such that the following assertions hold.
\begin{enumerate}
\item[(i)] If condition $(G1)$ is satisfied, then $\widehat{F}(x)\neq 0$ for $\|x\|_H \ge R_0$ and
\begin{equation*}
\mathrm{deg_B}(\widehat{F}, U_R^0) = 1 \qquad\mathrm{for}\quad R\ge R_0.
\end{equation*}
\item[(ii)] If condition $(G2)$ is satisfied, then $\widehat{F}(x)\neq 0$ for $\|x\|_H \ge R_0$ and
\begin{equation*}
\mathrm{deg_B}(\widehat{F}, U_R^0) = (-1)^{\dim X_0} \qquad\mathrm{for}\quad R\ge R_0.
\end{equation*}
\end{enumerate}
Here we denote $U_r^0 := \{x\in X_0 \ | \ \|x\|_H < r\}$ for $r > 0$.
\end{lemma}

\noindent\textbf{Proof of Lemma \ref{th-guid-fun-kk22}.} Let us note that from Lemma \ref{th-guid-fun-kk} we obtain the existence of $R_0 > 0$ with the property that
$\widehat{F}(x)\neq 0$ for $x\in X_0 $ with $\|x\|_H \ge R_0$ and
\begin{equation}
\begin{aligned}\label{amama}
& \mathrm{deg_B}(\widehat{F}, U_R^0) = 1 && \mathrm{for}\quad R\ge R_0 \text{ if condition $(G1)$ holds,} \\
& \mathrm{deg_B}(\widehat{F}, U_R^0) = (-1)^{\dim X_0} && \mathrm{for}\quad R\ge R_0 \text{ if condition $(G2)$ holds,}
\end{aligned}
\end{equation}
We show that $R_0$ is also the constant that we seek in lemma. To this end, assume that $R > R_0$. Given $\varepsilon > 0$, let us introduce the maps ${\bf A}_\varepsilon, \widehat{{\bf F}}_\varepsilon:{\bf E}_0\to{\bf E}_0$ given by
$${\bf A}_\varepsilon (x,y) := {\bf A}_0(x,y) + (0, \varepsilon x) \text{ \ and \ } \widehat{{\bf F}}_\varepsilon(x,y) := (0,\varepsilon x + \widehat{F}(x)) \text{ \ for \ } (x,y)\in{\bf E}_0.$$
One can easily check that $$\mathrm{Ker}\,(\mu I - {\bf A}_\varepsilon) = \{0\} \qquad \mathrm{for} \quad \mu\in (-\infty, 0]$$ and hence $(-\infty, 0]\subset\varrho({\bf A}_\varepsilon)$. Therefore we can define
$\widetilde H:[0,1]\times \overline{U_R} \to{\bf E}_0$ by $$\widetilde H(\mu, (x,y)) := \mu {\bf A}_\varepsilon (x,y) + (1 - \mu)(x,y) - (\mu I + (1 - \mu) {\bf A}_\varepsilon^{-1})\widehat{{\bf F}}_\varepsilon(x,y)$$ for $\mu\in [0,1]$, $(x,y)\in \overline{U_R}$. We claim that $\widetilde H$ has no zeros on $\partial U_R$. Suppose contrary to our claim that $\widetilde H(\mu,(x,y)) = 0$ for some $\mu\in [0,1]$ and $(x,y)\in\partial U_R$. If $\mu = 0$ then $(x,y) = {\bf A}_\varepsilon^{-1}\widehat{{\bf F}}_\varepsilon(x,y)$ and consequently $-{\bf A}_0 (x,y) + \widehat{{\bf F}}(x,y) = 0$. Hence $y = 0$ and $\widehat{F}(x) = 0$, where $\|x\|_H = R > R_0$ by \eqref{wz-norm}. This contradicts the choice of the number $R_0$. On the other hand, if $\mu\in(0,1]$ then
\begin{align*}
& (1/\mu - 1)(x,y) + {\bf A}_\varepsilon (x,y) = (I + (1/\mu - 1) {\bf A}_\varepsilon^{-1})\widehat{{\bf F}}_\varepsilon(x,y), \\
& (x,y) = ((1/\mu - 1)I + {\bf A}_\varepsilon)^{-1} (I + (1/\mu - 1) {\bf A}_\varepsilon^{-1})\widehat{{\bf F}}_\varepsilon(x,y),
\end{align*}
which, by the resolvent identity, gives $(x,y) = {\bf A}_\varepsilon^{-1} \widehat{{\bf F}}_\varepsilon(x,y)$. Similarly as before we deduce that $y = 0$ and $\widehat{F}(x) = 0$, where $\|x\|_H = R > R_0$. This is again a contradiction and our claim follows.
Using the claim and the homotopy invariance
\begin{equation}
\begin{aligned}\label{njnj4}
\mathrm{deg_B}({\bf A}_0 - \widehat{{\bf F}}, U_R) & = \mathrm{deg_B}({\bf A}_\varepsilon - \widehat{{\bf F}}_\varepsilon, U_R) = \mathrm{deg_B}(H(1,\,\cdot\,), U_R) \\
& = \mathrm{deg_B}(H(0,\,\cdot\,), U_R) = \mathrm{deg_B}(I - {\bf A}_\varepsilon^{-1}{\bf F}_\varepsilon, U_R).
\end{aligned}
\end{equation}
In view of \eqref{wz-norm} it follows that $U_R\subset U_R^0\times U_R^0$ and furthermore $${\bf A}_\varepsilon^{-1}{\bf F}_\varepsilon(x,y) = \widetilde H (0, (x,y))\neq (x,y) \qquad \mathrm{for} \quad (x,y)\in {\bf E}_0\setminus U_R.$$ Therefore, by the addition property of Brouwer degree
\begin{equation}\label{njnj3}
    \mathrm{deg_B}(I - {\bf A}_\varepsilon^{-1}{\bf F}_\varepsilon, U_R) = \mathrm{deg_B}(I - {\bf A}_\varepsilon^{-1}{\bf F}_\varepsilon, U_R^0\times U_R^0).
\end{equation}
On the other hand, let us observe that a simply calculations show that $$(I - {\bf A}_\varepsilon^{-1}{\bf F}_\varepsilon)(x,y) = (-1/\varepsilon\widehat{F}(x), y) \qquad \mathrm{for} \quad (x,y)\in{\bf E}_0.$$ Combining this with \eqref{njnj3}, \eqref{njnj4} and the multiplication property we infer that
\begin{equation*}
\begin{aligned}
   \mathrm{deg_B}({\bf A}_0 - \widehat{{\bf F}}, U_R) & = \mathrm{deg_B}(I - {\bf A}_\varepsilon^{-1}{\bf F}_\varepsilon, U_R) = \mathrm{deg_B}(I - {\bf A}_\varepsilon^{-1}{\bf F}_\varepsilon, U_R^0\times U_R^0)  \\
    & = \mathrm{deg_B}(-\widehat{F}, U_R^0) \cdot \mathrm{deg_B}(I, U_R^0) = (-1)^{\dim X_0} \, \mathrm{deg_B}(\widehat{F}, U_R^0),
\end{aligned}
\end{equation*}
which together with \eqref{amama} completes the proof of lemma. \hfill $\square$

\subsection{Proof of Theorem \ref{th-reso-m-h}.}
Let us first observe that from Lemma \ref{lem-bounkk} it follows that there is a constant $R_1 > 0$ such that for any $T$-periodic solution $w=(u,v):[0,+\infty)\to{\bf E}$ of the equation \eqref{diff:2kkaabb} one has
\begin{equation}\label{nier1111}
    \|{\bf Q} w(t)\|_{\bf E} \le R_1 \qquad \mathrm{for} \quad t\in [0,T],
\end{equation}
which, by \eqref{norm12} and \eqref{proj}, implies that
\begin{equation*}
    \|Q u(t)\|_\alpha \le R_1 \qquad \mathrm{for} \quad t\in [0,T].
\end{equation*}
Using Lemma \ref{lem-per-bound22}, we obtain a constant $R_2 > 0$ such that the equation \eqref{diff:2kkaabbcc} does not admit $T$-periodic solutions $w:[0,+\infty)\to {\bf E}$ such that $\|{\bf P} w(0)\|_{{\bf E}_0} \ge R_2$.
Furthermore, from Lemma \ref{th-guid-fun-kk22} it follows that there is $R_3 > R_2$ such that
$-{\bf A}_0 z + \widehat{{\bf F}}(z) \neq 0$ for $z\in {\bf E}_0 $ with $\|z\|_{{\bf E}_0} \ge R_3$ and furthermore
\begin{equation}
\begin{aligned}\label{roror1}
& \mathrm{deg_B}({\bf A}_0 - \widehat{{\bf F}}, U_R) = (-1)^{\dim X_0} && \mathrm{for}\quad R\ge R_3\text{ \ if $(G1)$ holds,}\\
& \mathrm{deg_B}({\bf A}_0 - \widehat{{\bf F}}, U_R) = 1 && \mathrm{for}\quad R\ge R_3 \text{ \ if $(G2)$ holds}.
\end{aligned}
\end{equation}
Let us define the following sets
\begin{align*}
U & :=\{(x,y)\in{\bf E}_0\ | \ \|(x,y)\|_{{\bf E}_0} < R_3 + 1\}, \\
V & :=\{(x,y)\in{\bf E}_-\oplus{\bf E}_+ \ | \ \|(x,y)\|_{\bf E} < R_1 + 1\},
\end{align*}
where we remind that the norm $\|\cdot\|_{{\bf E}_0}$ is given by \eqref{wz-norm} and $\|(x,y)\|_{\bf E}:= \|x\|_\alpha + \|y\|$ for $(x,y)\in {\bf E}$. Since the norms $\|\cdot\|_{\bf E}$ and $|\cdot|$ are equivalent and ${\bf E}_0$ is a finite dimensional space, it follows that the set $U\oplus V$ is open in the space ${\bf E}$, which is equipped with the norm $|\cdot|$.
We claim that
\begin{equation*}
{\bf \Psi}_T(s, (x,y)) \neq (x,y) \qquad \mathrm{for} \quad s\in[0,1], \ (x,y)\in\partial (U\oplus V).
\end{equation*}
Indeed, if the claim is false, then there is $s\in [0,1]$ and a mild solution $w:[0,+\infty)\to {\bf E}$ of the equation \eqref{diff:2kkaabb} such that $w(0) = w(T)\in\partial (U\oplus V)$.
Then either ${\bf P} w(0) = {\bf P} w(T)\in \partial U$ or ${\bf Q} w(0) = {\bf Q} w(T)\in \partial V$. Since the inequality \eqref{nier1111} holds, one has
\begin{equation*}
    \|{\bf P} w(0)\|_{{\bf E}_0} = R_3 + 1 \ge R_2,
\end{equation*}
which is impossible in view of the choice of the number $R_2$ and the claim follows. Therefore, by the homotopy invariance of topological degree
\begin{equation}
\begin{aligned}\label{jijiji}
\mathrm{deg_C}(I - {\bf \Phi}_T, U\oplus V) & = \mathrm{deg_C}(I - {\bf \Psi}_T(1,\,\cdot\,), U\oplus V) \\
& = \mathrm{deg_C}(I - {\bf \Psi}_T(0,\,\cdot\,), U\oplus V).
\end{aligned}
\end{equation}
In the remain part of the proof we will show that
\begin{equation}
\begin{aligned}\label{eq1122}
& \mathrm{deg_C}(I - {\bf \Psi}_T(0,\,\cdot\,), U\oplus V) = (-1)^{d_k} && \text{if $(G1)$ holds and} \\
& \mathrm{deg_C}(I - {\bf \Psi}_T(0,\,\cdot\,), U\oplus V) = (-1)^{d_{k-1}} && \text{if $(G2)$ holds.}
\end{aligned}
\end{equation}
Combining \eqref{jijiji} and \eqref{eq1122}, the proof of theorem will be completed if only we take $W:= U\oplus V$. To prove \eqref{eq1122} we will consider on the space ${\bf E}_0$, the following differential equation
\begin{equation}\label{bhbh}
\dot w(t) =  -{\bf A}_0 w(t) + {\bf P}{\bf F}(t, w(t)), \qquad  t > 0.
\end{equation}
If we denote by $\psi_T:{\bf E}_0\to{\bf E}_0$ the associated Poincar\'e operator, then using Corollary \ref{cor-dec-2aba} $(i)$, it is not difficult to see that $${\bf \Psi}_T(0, z) = S_{\bf A}(T)z_- + S_{\bf A}(T)z_+ + \psi_T(z_0) \qquad \mathrm{for} \quad z\in{\bf E},$$ where $z_\pm = {\bf Q}_\pm z$, $z_0 = {\bf P} z$. Let us denote by $H:[0,1]\times{\bf E}\to{\bf E}$ the map given by
$$H(\tau,z):= \tau S_{\bf A}(T)z_+ + S_{\bf A}(T)z_- + \psi_T(z_0) \qquad \mathrm{for} \quad (\tau, z)\in [0,1]\times {\bf E}.$$
The rest of the proof is divided into five steps. \\[5pt]
\textbf{Step 1.} We claim that $H$ is a condensing homotopy, that is,
$$\beta(H([0,1]\times \Omega)) \le e^{-\delta T} \beta(\Omega)$$
for any bounded $\Omega\subset{\bf E}$, where $\delta > 0$ is a constant. Indeed, if $\Omega\subset{\bf E}$ bounded, then using properties of Hausdorff measure of noncompactness (see Appendix), one has
\begin{align*}
\beta(H([0,1]\times \Omega)) & \le \beta(\{\tau S_{\bf A}(T){\bf Q}_+z \ | \ \tau\in[0,1], \ z\in\Omega\}) \\
& \qquad + \beta(S_{\bf A}(T){\bf Q}_-\Omega) + \beta(\psi_T({\bf P}\Omega)) \\
& = \beta(\{\tau S_{\bf A}(T){\bf Q}_+z \ | \ \tau\in[0,1], \ z\in\Omega\}),
\end{align*}
where the last inequality follows from the fact that the sets $S_{\bf A}(T){\bf Q}_-\Omega$ and $\psi_T({\bf P}\Omega)$ are relatively compact as continuous images of relatively compact sets. Since $$\{\tau S_{\bf A}(T){\bf Q}_+z \ | \ \tau\in[0,1], \ z\in\Omega\}\subset\mathrm{conv} ((S_{\bf A}(T){\bf Q}_+\Omega) \cup\{0\}),$$ by the properties of the Hausdorff measure of noncompactness again, we infer that
\begin{equation}
\begin{aligned}\label{rrrrr1}
& \beta(\{\tau S_{\bf A}(T){\bf Q}_+z \ | \ \tau\in[0,1], \ z\in\Omega\})  \le \beta(\mathrm{conv} ((S_{\bf A}(T){\bf Q}_+\Omega) \cup\{0\})) \\
& \hspace{37mm} = \beta((S_{\bf A}(T){\bf Q}_+\Omega) \cup\{0\}) = \beta(S_{\bf A}(T){\bf Q}_+\Omega).
\end{aligned}
\end{equation}
Therefore, by \eqref{rrrrr1}, \eqref{rownorm} and \eqref{rnon1ab}, we deduce that
\begin{equation*}
\begin{aligned}
\beta(H([0,1]\times \Omega)) \le \beta(S_{\bf A}(T){\bf Q}_+\Omega) \le e^{-\delta T}\beta({\bf Q}_+\Omega) \le e^{-\delta T}\beta(\Omega)
\end{aligned}
\end{equation*}
and the claim follows. \\[5pt]
\textbf{Step 2.} We prove that $H(\tau,z)\neq z$ for $\tau\in[0,1]$ and $z\in\partial (U\oplus V)$. To see this suppose that $H(\tau,z) = z$ for some $\tau\in[0,1]$ and $z\in\partial(U\oplus V)$. Then, from Corollary \ref{cor-dec-2aba} $(i)$ it follows that
\begin{align}\label{nhnh2}
\tau S_{\bf A}(T)z_+ + S_{\bf A}(T)z_- = z_+ + z_-  \quad\text{and}\quad \psi_T(z_0) = z_0.
\end{align}
Let us note that every mild solution of the equation \eqref{bhbh} is also a mild solution of the equation \eqref{diff:2kkaabbcc} with the parameters $\mu=1$ and $s = 0$. Using this fact together with the inequality $R_3 + 1 > R_2$, we infer that $\psi_T(z) \neq z$ for $z\in\partial U$ which implies that $z_- + z_+\in \partial V$ and furthermore
\begin{align*}
\tau S_{\bf A}(T)z_+ = z_+ \quad\text{and}\quad S_{\bf A}(T)z_- =  z_-
\end{align*}
The second equality together with Theorem \ref{th-spec} $(iii)$, give $z_-\in\mathrm{Ker}\,{\bf A} = \mathrm{Ker}\,(\lambda I - A)\times \{0\}$ and consequently $z_-\in{\bf E}_0$. But $z_-\in{\bf E}_-$, which implies that $z_- = 0$. Then $z_+ \in \partial V$ and $\tau S_{\bf A}(T)z_+ = z_+$. Iterating this equation $k$-times we obtain $\tau^k S_{\bf A}(kT)z_+ = z_+$ for $k\ge 1$. Hence, using \eqref{nieee1}, we obtain $$\|z_+\|_{\bf E} = \|\tau^k S_{\bf A}(kT)z_+\|_{\bf E} \le M \tau^k e^{-\delta k T} \|z_+\|_{\bf E} \qquad \mathrm{for} \quad k\ge 1.$$ This implies that $z_+ = 0$ and consequently $0 = z_+ + z_-\in \partial V$. This contradicts definition of $V$ and the claim follows. \\[5pt]
\textbf{Step 3.} In this step we prove that
\begin{equation}\label{degfor5}
\begin{aligned}
\mathrm{deg_C}(I - {\bf \Psi}_T(0,\,\cdot\,), U\oplus V) = (-1)^{d_{k-1}} \cdot \mathrm{deg_B}(I - \psi_T, U).
\end{aligned}
\end{equation}
To this end, let us note that by Step 1 and Step 2, $H$ is an admissible homotopy and consequently
\begin{equation}
\begin{aligned}\label{de1kkdeg}
\mathrm{deg_C}(I - {\bf \Psi}_T(0,\,\cdot\,), U\oplus V) & = \mathrm{deg_C}(I - H(1,\,\cdot\,), U\oplus V) \\
& = \mathrm{deg_C}(I - H(0,\,\cdot\,), U\oplus V).
\end{aligned}
\end{equation}
Since ${\bf E}_-$ is a finite dimensional space, the auxiliary linear operator $L:{\bf E}_+\oplus{\bf E}_-\to {\bf E}_+\oplus{\bf E}_-$ given by $L(z) = S_{\bf A}(T)z_-$, is compact and $\mathrm{Ker}\,(I - L) = \{0\}$, as we shown in Step 2. Hence, by the multiplication property of Leray-Schauder degree
\begin{equation}
\begin{aligned}\label{ro-dod}
\mathrm{deg_C}(I - H(0,\,\cdot\,), U\oplus V) & = \mathrm{deg_{LS}}(I - H(0,\,\cdot\,), U\oplus V) \\
& = \mathrm{deg_{LS}}(I - L, V) \cdot \mathrm{deg_B}(I - \psi_T, U).
\end{aligned}
\end{equation}
If $\lambda = \lambda_k$ with $k=1$, then by Theorem \ref{th-spec-dec} one has ${\bf E}_- = \{0\}$ and $L = 0$. Hence, by \eqref{ro-dod}, we obtain \eqref{degfor5}. If $k \ge 2$ then denoting
$$V_-:=\{z\in {\bf E}_-\ | \ \|z\|_{\bf E} < R_1 + 1\}, \quad V_+:=\{z\in {\bf E}_+\ | \ \|z\|_{\bf E} < R_1 + 1\}$$ and using addition and multiplication property of topological degree
\begin{equation}
\begin{aligned}\label{degform4}
\mathrm{deg_{LS}}(I - L, V) & = \mathrm{deg_{LS}}(I - L, V_-\oplus V_+) \\
& = \mathrm{deg_{LS}}(I, V_+)\cdot \mathrm{deg_B}(I - S_{{\bf A}}(T)|_{{\bf E}_-}, V_-) \\
& = \mathrm{deg_B}(I - S_{{\bf A}}(T)|_{{\bf E}_-}, V_-) = (-1)^{m_0},
\end{aligned}
\end{equation}
where $m_0$ is the sum of algebraic multiplicities of real eigenvalues of the operator $S_{{\bf A}}(T)|_{{\bf E}_-}$ which are greater than $1$. To determine the number $m_0$ precisely, let us note that, by Theorem \ref{th-spec-dec} $(i)$, one has ${\bf E}_- = K_1^-\oplus\ldots \oplus K_{k-1}^-$. Furthermore, by Lemma \ref{lem-22}, we have the inclusions $$K_i^-\subset \mathrm{Ker}\,(\mu^-_i I - {\bf A})\subset \mathrm{Ker}\,(e^{-\mu^-_i T} I - S_{\bf A}(T)) \qquad \mathrm{for} \quad i=1,2,\ldots, k-1,$$
where the numbers $\mu^-_i$ are from Theorem \ref{th-spec}.
Using this with Lemma \ref{lem-wart}, gives $$\sigma(S_{{\bf A}}(T)|_{{\bf E}_-}) = \{e^{-\mu^-_i T} \ | \ 1\le i \le k-1\}$$ and the algebraic multiplicity of each eigenvalue $e^{-\mu^-_i T}$ is equal to $\dim K_i^-$. Since, $\mu^-_i < 0$ for $i=1,\ldots,k-1$ one has $$m_0 = \sum_{i=1}^{k-1} \dim K_i^- = \dim{\bf E}_-,$$ which together with \eqref{degform4} gives
\begin{equation}
\begin{aligned}\label{degform4a}
\mathrm{deg_{LS}}(I - L, V) = (-1)^{d_{k-1}}.
\end{aligned}
\end{equation}
Combining \eqref{de1kkdeg}, \eqref{ro-dod} and \eqref{degform4a} yields \eqref{degfor5} as desired.  \\[5pt]
\textbf{Step 4.} Let us prove that
\begin{equation}
\begin{aligned}\label{wz-ind}
    \mathrm{deg_B}(I - \psi_T, U) = \mathrm{deg_B}({\bf A}_0 - \widehat{{\bf F}}, U).
\end{aligned}
\end{equation}
For this purpose, consider the family differential equations of the form
\begin{equation}\label{bhbh2}
\dot w(t) = - \mu{\bf A}_0 w(t) + \mu{\bf P}{\bf F}(t, w(t)), \qquad  t > 0.
\end{equation}
where $\mu\in(0,1]$ is a parameter and let ${\bf \Theta}^\mu_T: {\bf E}_0\to{\bf E}_0$ be the Poincar\'e operator associated with this equation. Let us note that every mild solution of the equation \eqref{bhbh2} is also a mild solution of the equation \eqref{diff:2kkaabbcc} with the parameter $s=0$. From this fact and the definition of $U$ we deduce that ${\bf \Theta}^\mu_T(x,y)\neq (x,y)$ for $\mu\in(0,1]$ and $(x,y)\in \partial U$. Therefore, for any $\mu\in(0,1]$, we have
\begin{equation}
\begin{aligned}\label{degdeg}
    \mathrm{deg_B}(I - \psi_T, U) = \mathrm{deg_B}(I - {\bf \Theta}^1_T, U) = \mathrm{deg_B}(I - {\bf \Theta}^\mu_T, U).
\end{aligned}
\end{equation}
Observe that the neighborhood $U$ was chosen such that
\begin{equation*}
    -{\bf A}_0 (x,y) + \widehat{{\bf F}}(x,y) \neq 0 \qquad \mathrm{for} \quad (x,y)\in {\bf E}_0\setminus U.
\end{equation*}
Hence, by Theorem \ref{th-kras-2}, there is $\mu_0 \in (0,1)$ such that for any $\mu\in(0,\mu_0]$ one has
\begin{equation}\label{degfor1}
    \mathrm{deg_B}(I - {\bf \Theta}^\mu_T, U) = \mathrm{deg_B}({\bf A}_0 - \widehat{{\bf F}}, U).
\end{equation}
Therefore, combining \eqref{degdeg} with \eqref{degfor1} gives \eqref{wz-ind}. \\[5pt]
\textbf{Step 5.} Finally we are in the position to prove \eqref{eq1122}. To this end it is enough to combine \eqref{roror1} with the assertions from Steps 3 and 4. Thus the proof of the theorem is completed. \hfill $\square$ \\

\section{Applications}

In this section we provide applications of the abstract results obtained in the previous sections to partial differential equations. We will assume that $\Omega\subset\mathbb{R}^n$, $n\ge 1$, is an open bounded set with the boundary $\partial\Omega$ of class $C^1$. Consider the strongly damped wave equation
\begin{equation}\label{A-eps-res-ah}
\left\{\begin{aligned}
& u_{tt}= - c\mathcal{A} \, u_t - \mathcal{A} \, u + \lambda  u + f(t,x, u), && t > 0, \ x\in\Omega, \\
& u(t,x) = 0 && t \ge 0, \ x\in\partial\Omega,
\end{aligned}\right.
\end{equation}
where $c> 0$ is a damping factor, $\lambda$ is a real number and $\mathcal{A}$ is a differential operator of the following form
$$\mathcal{A} \bar u (x) = -\sum_{i,j=1}^n D_j(a_{ij}(x)D_i \bar u(x)) \qquad \mathrm{for} \quad \bar u\in C^1(\overline{\Omega}),$$ which is {\em symmetric} $a_{ij} = a_{ji}\in C^1(\overline\Omega)$ and {\em uniformly elliptic}, that is,
$$\sum_{1\le i,j\le n}a_{ij}(x)\xi^i\xi^j \ge c_0 |\xi|^2 \qquad \mathrm{for} \quad x\in\Omega, \ \xi\in\mathbb{R}^n, \ \ \text{where} \ \ c_0 > 0.$$
Furthermore assume that $f:[0,+\infty)\times\Omega\times\mathbb{R}\to\mathbb{R}$ is a continuous map such that \\[2pt]
\noindent\makebox[22pt][l]{$(E1)$} \parbox[t][][t]{118mm}{there is $L> 0$ such that if $t\in [0,+\infty)$,  $x\in\Omega$ and $s_1.s_2\in\mathbb{R}$, then
    \begin{align*}
    |f(t,x,s_1) - f(t,x,s_2)| & \le L |s_1 - s_2|;
    \end{align*}}\\
\noindent\makebox[22pt][l]{$(E2)$} \parbox[t][][t]{118mm}{there is $m > 0$ such that $$|f(t,x,s)| \le m \qquad \mathrm{for} \quad t\in [0,+\infty), \ x\in\Omega, \ s\in\mathbb{R};$$}\\[5pt]
\noindent\makebox[22pt][l]{$(E3)$} \parbox[t][][t]{118mm}{there is $T > 0$ such that $$f(t,x,s) = f(t + T,x,s) \ \ \text{for} \ t\in[0,+\infty), \ x\in\Omega, \ s\in\mathbb{R}.$$} \\
Let us introduce the abstract framework for the equation \eqref{A-eps-res-ah}. To this end, let us denote $X:=L^p(\Omega)$, for $p\ge 2$, and define the operator $A_p: X\supset D(A_p)\to X$ by
\begin{equation*}
\begin{aligned}
D(A_p) := W^{2,p}(\Omega) \cap W^{1,p}_0(\Omega), \quad A_p \bar u := \mathcal{A} \bar u \quad \text{for} \ \ \bar u\in D(A_p).
\end{aligned}
\end{equation*}
It is known (see e.g. \cite{MR1778284}, \cite{Pazy}, \cite{MR500580}) that $A_p$ is a positive sectorial operator. Let us denote by $X^\alpha := D(A_p^\alpha)$, for $\alpha\in(0,1)$, the associated fractional space and define the map $F\colon [0,+\infty)\times X^\alpha \to X$, given for any $\bar u\in X^\alpha$, by
\begin{equation*}
    F(t,\bar u)(x) := f(t,x, \bar u(x)) \qquad \mathrm{for} \quad t\in [0,+\infty), \ x\in\Omega.
\end{equation*}
We call $F$ \emph{the Nemitskii operator} associated with $f$. We are ready to write the equation \eqref{A-eps-res-ah} in the following abstract form
\begin{equation}\label{row-abs}
\ddot u(t)  = - A_p u(t) - cA_p \dot u(t) + \lambda u(t) + F (t,u(t)),  \qquad   t > 0 \\
\end{equation}

\begin{remark}\label{rem-pom2}
$(a)$ We claim that assumptions $(A1)-(A3)$ are satisfied. \\[2pt]
Indeed, $(A1)$ holds because $A_p$ has compact resolvent as it was proved for example in \cite{MR1778284}, \cite{Pazy}, \cite{MR500580}.
To see that $(A2)$ holds it is enough to take $H:=L^2(\Omega)$ equipped with the standard inner product and norm.
Since $p\ge 2$ and $\Omega$ is a bounded set, the embedding $i:L^p(\Omega) \hookrightarrow L^2(\Omega)$ is well-defined and continuous. Furthermore, if we define $\widehat A:= A_2$, then $$i(D(A_p))\subset D(\widehat A) \quad\text{and}\quad i(A_p \bar u) = \widehat A i(\bar u) \quad\text{for}\quad \bar u \in D(A_p),$$
which shows that $A_p \subset \widehat A$ in the sense of the inclusion $i\times i$. Since the operator $\widehat A$ is self-adjoint (see e.g. \cite{MR1778284}) we see that the assumption $(A3)$ is also satisfied. \\[5pt]
$(b)$ Let us observe that $F$ is satisfies assumptions $(F1)$, $(F3)$, $(F4)$ and $(F5)$. \\[2pt]
Since $f$ satisfies assumptions $(E1)-(E3)$, the fact that conditions $(F1)$, $(F4)$ and $(F5)$ hold is straightforward. We only show assumption $(F3)$. To this end take a bounded sequence $(\bar u_n)$ in $X^\alpha$. Since $A_p$ has compact resolvents, by \cite[Theorem 1.4.8]{MR610244}, the inclusion $X^\alpha\hookrightarrow X$ is compact, and hence, passing if necessary to a subsequence, we can assume that $\bar u_n \to \bar u_0$ in $X$ as $n\to +\infty$. Therefore, using $(E2)$ and the dominated convergence theorem, it is not difficult to verify that $F(\bar u_n) \to F(\bar u_0)$ in $X$ as $n\to +\infty$, which proves that $(F3)$ holds. \hfill $\square$
\end{remark}

\subsection{Properties of the Nemitskii operator}

We proceed to examine when the the Nemitskii operator $F$ satisfies geometrical conditions $(G1)$ and $(G2)$. Let us first note that, by Remark \ref{rem-pom}, the spectrum $\sigma(A_p)$ of the operator $A_p$ consists of the sequence of positive eigenvalues $$0 < \lambda_1 < \lambda_2 < \ldots < \lambda_i < \lambda_{i+1} < \ldots$$ which is finite or $\lambda_i \to +\infty$ as $i\to +\infty$. We recall also that $X^\alpha_+$, $X^\alpha_-$ and $X_0$ are subspaces obtained in Remark \ref{rem-proj}, but this time for the operator $A_p$. In particular $X_0 = \mathrm{Ker}\,(\lambda I - A_p)$. Let us start with the following theorem which says that the conditions $(G1)$ and $(G2)$ are implicated by the well-known \emph{Landesman-Lazer} conditions introduced in \cite{MR0267269}.

\begin{theorem}\label{lem-est2}
Let $f_+,f_-\colon [0,+\infty)\times\Omega \to \mathbb{R}$ be continuous functions such that
\begin{align*}
f_\pm(t,x) = \lim_{s \to \pm\infty} f(t,x,s) \quad\text{for}\quad x\in\Omega, \ \  \text{uniformly for} \ \ t\in[0,+\infty).
\end{align*}
\makebox[5mm][l]{(i)} \parbox[t]{117mm}{Assume that for every $t\in[0,T]$ the following holds
\begin{equation*}
\int_{\{u>0\}} f_+(t,x) \bar u(x) \,d x  + \int_{\{u<0\}} f_-(t,x) \bar u(x) \,d x > 0 \ \ \text{for} \ \ \bar u\in X_0\setminus\{0\}.\leqno{(LL1)}
\end{equation*}
If the sets $B_1\subset X^\alpha_+\oplus X^\alpha_-$ and $B_2 \subset X_0$ are bounded in the norms $\|\cdot\|_\alpha$ and $\|\cdot\|_{L^2}$, respectively, then there is a constant $R> 0$ such that
\begin{equation*}
\<F(t,\bar w + \bar u), \bar u\>_{L^2} > - \<F(t,\bar w + \bar u), \bar v\>_{L^2}
\end{equation*}
for any $t\in[0,T]$ and $(\bar w, \bar v,\bar u)\in B_1\times B_2 \times X_0$, with $\|\bar u\|_{L^2} \ge R$.}\\[5pt]
\makebox[5mm][l]{(ii)} \parbox[t]{117mm}{Assume that for every $t\in[0,T]$ the following holds
\begin{equation*}
\int_{\{u>0\}} f_+(t,x) \bar u(x) \,d x  + \int_{\{u<0\}} f_-(t,x) \bar u(x) \,d x < 0 \ \ \text{for} \ \ \bar u\in X_0\setminus\{0\}.\leqno{(LL2)}
\end{equation*}
If the sets $B_1\subset X^\alpha_+\oplus X^\alpha_-$ and $B_2 \subset X_0$ are bounded in the norms $\|\cdot\|_\alpha$ and $\|\cdot\|_{L^2}$, respectively, then there is a constant $R> 0$ such that
\begin{equation*}
\<F(t,\bar w + \bar u), \bar u\>_{L^2} < - \<F(t,\bar w + \bar u), \bar v\>_{L^2}
\end{equation*}
for any $t\in[0,T]$ and $(\bar w, \bar v,\bar u)\in B_1\times B_2 \times X_0$, with $\|\bar u\|_{L^2} \ge R$.}
\end{theorem}
\noindent\textbf{Proof.} Except for technical modifications, the argument goes in the lines of the proof of
\cite[Theorem 4.3]{Kok5}. We encourage the reader to reconstruct details. \hfill $\square$ \\[-5pt]

The following theorem shows that the conditions $(G1)$ and $(G2)$ are also implicated by \emph{the strong resonance conditions}, studied for example in \cite{MR713209}, \cite{MR597281}.

\begin{theorem}\label{lem-est3}
Let $f_\infty \colon [0,+\infty)\times \Omega \to \mathbb{R}$ be a continuous function such that
\begin{equation}\label{asu}
f_\infty(t,x)  = \lim_{|s| \to +\infty} f(t,x,s)\cdot s \ \ \text{for}\ \  x\in\Omega, \ \ \text{uniformly for} \ \ t\in[0,+\infty).
\end{equation}
\makebox[6mm][l]{(i)}\parbox[t]{117mm}{Assume that \\[-9pt]
\begin{equation*}\leqno{(SR1)}
\quad\left\{\begin{aligned}
\parbox[t][][t]{102mm}{\em there is a function $h\in L^1(\Omega)$ such that \\[3pt] $f(t,x,s)\cdot s \ge h(x)$ for $(t,x,s)\in [0,+\infty)\times\Omega\times\mathbb{R}$ and \\[3pt]  $\int_\Omega f_\infty(t,x)\, d x > 0$ for $t\in[0,T]$.}
\end{aligned} \right.
\end{equation*}\\
If the sets $B_1\subset X^\alpha_+\oplus X^\alpha_-$ and $B_2 \subset X_0$ are bounded in the norms $\|\cdot\|_\alpha$ and $\|\cdot\|_{L^2}$, respectively, then there is a constant $R > 0$ such that
$$\<F(t,\bar w + \bar u), \bar u\>_{L^2} > - \<F(t,\bar w + \bar u), \bar v\>_{L^2}$$
for any $t\in[0,T]$ and $(\bar v, \bar w,\bar u)\in B_1\times B_2 \times X_0$, with $\|\bar u\|_{L^2}\ge R$.}\\[5pt]
\makebox[6mm][l]{(ii)}\parbox[t]{117mm}{Assume that \\[-9pt]
\begin{equation*}\leqno{(SR2)}
\quad\left\{\begin{aligned}
\parbox[t][][t]{102mm}{\em there is a function $h\in L^1(\Omega)$ such that \\[3pt] $f(t,x,s)\cdot s \le h(x)$ for $(t,x,s)\in [0,+\infty)\times\Omega\times\mathbb{R}$ and \\[3pt]  $\int_\Omega f_\infty(t,x)\, d x < 0$ for $t\in[0,T]$.}
\end{aligned} \right.
\end{equation*}\\
If the sets $B_1\subset X^\alpha_+\oplus X^\alpha_-$ and $B_2 \subset X_0$ are bounded in the norms $\|\cdot\|_\alpha$ and $\|\cdot\|_{L^2}$, respectively, then there is constant $R> 0$ such that
$$\<F(t,\bar w + \bar u), \bar u\>_{L^2} < - \<F(t,\bar w + \bar u), \bar v\>_{L^2}$$
for any $t\in[0,T]$ and $(\bar v, \bar w,\bar u)\in B_1\times B_2 \times X_0$, with $\|\bar u\|_{L^2}\ge R$.}
\end{theorem}
\noindent\textbf{Proof.} Similarly as before, with some technical modifications, the argument goes in the lines of the proof of \cite[Theorem 4.5]{Kok5}. We encourage the reader one more time to reconstruct details. \hfill $\square$ \\

\begin{remark}
Let us observe that under the assumptions of Theorem \ref{lem-est3} one has
\begin{equation*}
f_\pm (t,x) = \lim_{s \to \pm\infty} f(t,x,s) = 0 \text{ \ for \ } x\in\Omega,
\end{equation*}
which implies that the Landesman-Lazer conditions $(LL1)$ and $(LL2)$ used in Theorem \ref{lem-est2} are not valid in this case. \hfill $\square$
\end{remark}

\subsection{Criteria on existence of periodic solutions}

Let us observe that the equation \eqref{row-abs} can by written as
\begin{equation}\label{row-dr2}
\dot w(t) = -{\bf A}_p w(t) + {\bf F} (t,w(t)),  \qquad   t > 0.
\end{equation}
where ${\bf A}_p:{\bf E}\supset D({\bf A}_p)\to{\bf E}$ is a linear operator on ${\bf E}:= X^\alpha\times X$ given by
\begin{equation*}
\begin{aligned}
D({\bf A}_p) & :=\{(\bar u,\bar v)\in X^\alpha\times X \ | \ \bar u + c \bar v\in D(A_p)\} \\
{\bf A}_p(\bar u,\bar v) & :=(-\bar v, A_p(\bar u + c \bar v) - \lambda \bar u)
\end{aligned}
\end{equation*}
and ${\bf F}: [0,+\infty)\times{\bf E}\to{\bf E}$ is a map defined by ${\bf F}(t,(\bar u,\bar v)):= (0, F(t,\bar u))$ for $t\in[0,+\infty)$ and $(\bar u,\bar v)\in {\bf E}$. Similarly as before we assume that the space ${\bf E}$ is equipped with the norm $|\cdot|$, obtained in Proposition \ref{lem-con-sem}. Then Remark \ref{rem-pom2} $(b)$ and Proposition \ref{th-exist1} $(a)$ assert that \emph{the Poincar\'e operator} ${\bf \Phi}_T: {\bf E} \to {\bf E}$ associated with the equation \eqref{row-dr2} is well-defined, continuous and
\begin{align*}
\beta({\bf \Phi}_T(\Omega)) \le e^{-\delta T}\beta(\Omega) \ \ \text{for any bounded} \ \ \Omega\subset{\bf E},
\end{align*}
where $\delta > 0$ is a constant. Let us first prove the following {\em criterion with Landesman-Lazer type conditions}.
\begin{theorem}\label{th-crit-ogrh}
Let $f_+,f_-\colon [0,+\infty)\times\Omega \to \mathbb{R}$ be continuous functions such that
\begin{equation*}
f_\pm(t,x) = \lim_{s \to \pm\infty} f(t,x,s) \quad\text{for}\quad x\in\Omega, \ \  \text{uniformly for} \ \ t\in[0,+\infty).
\end{equation*}
If $\lambda=\lambda_k$ is the $k$-th eigenvalue of the operator $A_p$, then there is an open set $W \subset {\bf E}$ such that ${\bf \Phi}_T(x,y)\neq (x,y)$ for $(x,y)\in\partial W$ and: \\[-10pt]
\begin{enumerate}
\item[(i)] $\mathrm{deg_{C}}(I - {\bf \Phi}_T, W) = (-1)^{d_k}$ \ if condition $(LL1)$ is satisfied, \\[-9pt]
\item[(ii)] $\mathrm{deg_{C}}(I - {\bf \Phi}_T, W) = (-1)^{d_{k-1}}$ \ if condition $(LL2)$ is satisfied,
\end{enumerate}
Here $d_l$ is such that $d_0 := 0$ and $d_l:= \sum_{i=1}^l \dim\mathrm{Ker}\,(\lambda_i I - A)$ for $l\ge 1$.
\end{theorem}
\noindent\textbf{Proof.} The proof is a consequence of Theorem \ref{th-reso-m-h} and Theorem \ref{lem-est2} \hfill $\square$ \\

Now we proceed to the following \emph{criterion with strong resonance conditions}. Similarly as in Theorem \ref{lem-est3} we make the restriction that $\Omega\subset\mathbb{R}^n$ where $n\ge 3$.
\begin{theorem}\label{th-crit-ogrsrh}
Let $f_\infty \colon [0,+\infty)\times\Omega \to \mathbb{R}$ be a continuous function such that
\begin{equation}\label{asu}
f_\infty(t,x)  = \lim_{|s| \to +\infty} f(t,x,s)\cdot s \ \ \text{for}\ \  x\in\Omega, \ \ \text{uniformly for} \ \ t\in[0,+\infty).
\end{equation}
If $\lambda=\lambda_k$ is the $k$-th eigenvalue of the operator $A_p$, then there is a neighborhood $W \subset {\bf E}$ such that ${\bf \Phi}_T(x,y)\neq (x,y)$ for $(x,y)\in\partial W$ and
\begin{enumerate}
\item[(i)] $\mathrm{deg_{C}}(I - {\bf \Phi}_T, W) = (-1)^{d_k}$ \ if condition $(SR1)$ is satisfied, \\[-9pt]
\item[(ii)] $\mathrm{deg_{C}}(I - {\bf \Phi}_T, W) = (-1)^{d_{k-1}}$ \ if condition $(SR2)$ is satisfied.
\end{enumerate}
Here $d_0 := 0$ and $d_l:= \sum_{i=1}^l \dim\mathrm{Ker}\,(\lambda_i I - A)$ for $l\ge 1$.
\end{theorem}
\noindent\textbf{Proof.} The proof is a consequence of Theorem \ref{th-reso-m-h} and Theorem \ref{lem-est3} \hfill $\square$

\begin{remark}
By Theorems \ref{th-crit-ogrh} and \ref{th-crit-ogrsrh} and existence property of topological degree it follows that the equation \eqref{row-dr2} admits a $T$-periodic mild solution provided either Landesman-Lazer or strong resonance conditions are satisfied.
\end{remark}

\section{Appendix}

\subsection{The Brouwer degree}
Consider the finite dimensional space $X$ and let $U\subset X$ be an open bounded set. For a continuous map $f:\overline U\to X$, such that $f(x)\neq 0$ for $x\in\partial U$, one can assign the integer number $\mathrm{deg_B}(f,U)$, called \emph{the Brouwer degree}, with the following properties. \\[5pt]
\makebox[8mm][l]{(B1)} \parbox[t]{118mm}{(Existence) If $\mathrm{deg_B}(f, U)\neq 0$ then there is
    $x\in U$ such that $f(x) = 0$.}\\[5pt]
\makebox[8mm][l]{(B2)} \parbox[t]{118mm}{(Additivity) If $f\colon\overline{U} \to X$ is such that
     $f(x)\neq 0$ for $x\in\partial U$ and $U_1, U_2\subset U$ are open disjoint sets such that $\{x\in\overline{U} \ | \ f(x) = 0\}\subset U_1\cup U_2$, then
    $$\mathrm{deg_B}(f, U) = \mathrm{deg_B}(f_{|\overline{U}_1}, U_1) + \mathrm{deg_B}(f_{|\overline{U}_2}, U_2).$$}\\
\makebox[8mm][l]{(B3)} \parbox[t]{118mm}{(Homotopy invariance) If a continuous map $h: [0,1]\times\overline U\to X$ is such that $h(\lambda,x)\neq 0$ for $(\lambda,x)\in[0,1]\times\partial U$, then
$$\mathrm{deg_B}(h(0,\,\cdot\,), U) = \mathrm{deg_B}(h(1,\,\cdot\,), U).$$}\\
\makebox[8mm][l]{(B4)} \parbox[t]{118mm}{(Normalization) If $0\in U$ then $\mathrm{deg_B}(I, U) = 1$.}\\[5pt]
\makebox[8mm][l]{(B5)} \parbox[t]{118mm}{(Multiplication) Let $U$ and $V$ be open bounded subsets of finite dimensional spaces $X$ and $Y$, respectively and let $f:\overline{U} \to X$ and $g:\overline{V} \to Y$ be continuous maps such that $f(x)\neq 0$ for $x\in\partial U$ and $g(y)\neq 0$ for $y\in\partial V$. Then $(f\times g)(x,y) \neq 0$ for $(x,y)\in\partial (U\times V)$ and $$\mathrm{deg_B}(f\times g, U\times V) = \mathrm{deg_B}(f, U)\cdot \mathrm{deg_B}(g, V).$$}\\[5pt]

\noindent Let us consider the following family of differential equations
\begin{equation}\label{rrrr}
\dot u(t) = \mu f(t,u(t)), \qquad t > 0
\end{equation}
where $\mu\in[0,1]$ is a parameter and $f:[0,+\infty)\times X \to X$ is a continuous bounded map. It is well-known that, for every initial data $x\in X$ and parameter $\mu\in(0,1]$, there is a smooth solution $u(\,\cdot\,;\mu,x):[0,+\infty)\to X$ of the equation \eqref{rrrr}. Given $T>0$, let $\varphi^\mu_T:\mathbb{R}^n\to\mathbb{R}^n$ be a Poincar\'e operator associated with this equation: $$\varphi^\mu_T(x) = u(T;\mu,x)  \ \ \text{for}  \ \ \mu\in(0,1] \ \ \text{and} \ \  x\in X.$$ The following result is a classical degree formula which connects the Brouwer degree of the operator $\varphi^\mu_T$ with the averaging of the right hand side of the equation \eqref{rrrr}.

\begin{theorem}{\em (see \cite{MR2171122}, \cite{Krasnoselskii})}\label{th-kras-2}
Assume that $U\subset X$ is an open bounded set such that $f(x)\neq 0$ for $x\in\partial U$. Then there is $\mu_0 > 0$ such that, if $\mu\in(0,\mu_0]$ then $\varphi^\mu_T(x) \neq x$ for $x\in\partial U$ and
\begin{equation*}
    \mathrm{deg_B}(I - \varphi^\mu_T, U) = \mathrm{deg_B}(-\widehat{f}, U),
\end{equation*}
where $\widehat{f}(x) := \frac{1}{T}\int_0^T f(\tau, x) \,d\tau$ for $x\in\mathbb{R}^n$ is an averaging of the map $f$.
\end{theorem}

\subsection{Hausdorff measure of noncompactness}\label{mesu}

We briefly recall the notion of Hausdorff measure of noncompactness. For more details see \cite{MR1189795}, \cite{Dugundji-Granas}. Let $X_0\subset X$ be a linear subspace of an infinite dimensional Banach space $X$ equipped with a norm $\|\cdot\|$. \emph{The Hausdorff measure of noncompactness} $\beta_{X_0}$ of bounded subset $\Omega\subset X_0$ is defined as
\begin{equation*}
\beta_{X_0}(\Omega):= \inf \{r > 0 \ | \ \Omega \subset \bigcup_{i=1}^{k_r} B(x_i, r), \ \text{ where } \ x_i\in X_0 \ \text{ for } \ i=1,\ldots, k_r\},
\end{equation*}
where $B(x, r) = \{y\in X\ | \ \|y - x\| < r\}$. The measure has the following properties: \\[-8pt]
\begin{itemize}
\item[(M1)] $\beta_{X_0}(\Omega) = 0$ if and only if $\Omega\subset X_0$ is a relatively compact set, \\[-8pt]
\item[(M2)] if $\Omega_1,\Omega_2$ are bounded sets and $\Omega_1 \subset \Omega_2$, then $\beta_{X_0}(\Omega_1) \le \beta_{X_0}(\Omega_2)$, \\[-8pt]
\item[(M3)] $\beta_{X_0}(\mathrm{conv}\,\Omega) = \beta_{X_0}(\Omega)$ for any bounded $\Omega\subset X_0$, \\[-8pt]
\item[(M4)] if $\Omega,\Omega_1,\Omega_2$ are bounded sets and $\lambda\in\mathbb{R}$, then
$$\qquad\beta_{X_0}(\lambda\Omega) = |\lambda| \beta_{X_0}(\Omega) \text{ \ and \ }
\beta_{X_0}(\Omega_1 + \Omega_2) \le \beta_{X_0}(\Omega_1) + \beta_{X_0}(\Omega_2).$$
\end{itemize}

If $X_0=X$ then for abbreviation we write $\beta:=\beta_{X_0}$. In the general situation the measures $\beta$ and $\beta_{X_0}$ are not equal. However we have the following lemma.
\begin{lemma}\label{lem-noncom}
Assume that there is a bounded linear map $P: X\to X$ with $P(X) = X_0$, $Px = x$ for $x\in X_0$ and $\|P\|\le 1$. Then, for any bounded set $\Omega\subset X_0$, we have $\beta_X(\Omega) = \beta_{X_0}(\Omega)$.
\end{lemma}
\noindent\textbf{Proof.} According to the definition of $\beta$, the inequality $\beta_{X_0}(\Omega)\ge\beta_X(\Omega)$ holds for any bounded set $\Omega\subset X_0$. To verify the opposite inequality, let us take $\varepsilon > 0$ and consider the covering of $\Omega$ by a finite number of balls $B(x_1, r_\varepsilon)$, $B(x_2, r_\varepsilon)$, \ldots, $B(x_n, r_\varepsilon)$ with radius $r_\varepsilon:=\beta(\Omega) + \varepsilon$. Then, for any $x\in B(x_i, r_\varepsilon)$, we have $\|Px - P x_i\| \le \|x - x_i\| \le r_\varepsilon$, which implies that $P B(x_i, r_\varepsilon) \subset B(P x_i, r_\varepsilon)$.
Consequently, for any $\Omega\subset X_0$, $$\Omega = P\Omega \subset P \bigcup_{i=1}^k B(x_i, r_\varepsilon) = \bigcup_{i=1}^k PB(x_i, r_\varepsilon) \subset \bigcup_{i=1}^k B(P x_i, r_\varepsilon)$$ and therefore the balls $B(P x_1, r_\varepsilon)$, $B(P x_2, r_\varepsilon)$, \ldots, $B(P x_n, r_\varepsilon)$ make a covering of $\Omega$ in $X_0$. Since $\varepsilon > 0$ is arbitrary small, it follows that $\beta_{X_0}(\Omega) \le \beta_X(\Omega)$ and the proof is completed. \hfill $\square$\\

\subsection{Topological degree for condensing vector fields}

We briefly describe the topological degree for condensing vector fields. For more details and construction see \cite{MR0246285}, \cite{MR0306986}, \cite{MR0312341}. Let $U\subset X$ be an open bounded subset of infinite dimensional Banach space $X$ and let $f\colon \overline U \to X$ be a continuous map. Given $k\in[0,1)$, we say that the map $f$ is \emph{$k$-condensing} provided
\begin{equation*}
\beta(f(\Omega))\le k\beta(\Omega) \quad \text{for any bounded set}  \ \ \Omega\subset \overline U,
\end{equation*}
where $\beta$ is the Hausdorff measure of noncompactness on $X$.
The map $h: [0,1]\times\overline U\to X$ is called \emph{a $k$-condensing homotopy} provided
\begin{equation*}
    \beta(h([0,1]\times\Omega))\le k\beta(\Omega) \quad \text{for any bounded set}  \ \ \Omega\subset \overline U.
\end{equation*}
We say that the vector field $I - f\colon\overline{U} \to X$ is admissible if $f\colon\overline{U} \to X$ is $k$-condensing and $f(x)\neq x$ for $x\in \partial U$.
For an admissible map $I - f:\overline{U} \to X$, we can assign the integer number $\mathrm{deg_C}(I - f, U)$, called {\em the topological degree for condensing vector fields}, satisfying the following properties. \\[3pt]
\makebox[8mm][l]{(C1)} \parbox[t]{118mm}{(Existence) If $\mathrm{deg_C}(I - f, U)\neq 0$, then there is
    $x\in U$ such that $F(x) = x$.}\\[5pt]
\makebox[8mm][l]{(C2)} \parbox[t]{118mm}{(Additivity) If $I - f\colon\overline{U} \to X$ is an admissible map and $U_1, U_2\subset U$ are disjoint open sets such that $\{x\in\overline{U} \ | \ F(x) = x\}\subset U_1\cup U_2$, then
    $$\mathrm{deg_C}(I - f, U) = \mathrm{deg_C}(I - f_{|\overline{U}_1}, U_1) + \mathrm{deg_C}(I - f_{|\overline{U}_2}, U_2).$$}\\
\makebox[8mm][l]{(C3)} \parbox[t]{118mm}{(Homotopy invariance) If $h: [0,1]\times\overline U\to X$ is a $k$-condensing homotopy such that $h(\lambda,x)\neq x$ for $(\lambda,x)\in[0,1]\times\partial U$, then
$$\mathrm{deg_C}(I - h(0,\,\cdot\,), U) = \mathrm{deg_C}(I - h(1,\,\cdot\,), U).$$}\\
\makebox[8mm][l]{(C4)} \parbox[t]{118mm}{(Normalization) If $0\in U$ then $\mathrm{deg_C}(I, U) = 1$.}

\def\cprime{$'$} \def\polhk#1{\setbox0=\hbox{#1}{\ooalign{\hidewidth
  \lower1.5ex\hbox{`}\hidewidth\crcr\unhbox0}}} \def\cprime{$'$}
  \def\cprime{$'$} \def\cprime{$'$}
\providecommand{\bysame}{\leavevmode\hbox to3em{\hrulefill}\thinspace}
\providecommand{\MR}{\relax\ifhmode\unskip\space\fi MR }
\providecommand{\MRhref}[2]{%
  \href{http://www.ams.org/mathscinet-getitem?mr=#1}{#2}
}
\providecommand{\href}[2]{#2}

\parindent = 0 pt

\end{document}